\documentclass[12pt]{article}
\usepackage{amssymb}
\usepackage{latexsym}
\usepackage{amsthm}
\usepackage{amscd}
\usepackage{amsmath}
\usepackage{diagrams}

\newtheorem{thm}{Theorem}[section]
\newtheorem{lem}[thm]{Lemma}

\newtheorem{defn-lem}[thm]{Definition-Lemma}

\theoremstyle{remark}
\newtheorem{rem}[thm]{Remark}
\newtheorem{claim}[thm]{Claim}
\theoremstyle{definition}
\newtheorem{defn}[thm]{Definition}
\numberwithin{equation}{section}

\allowdisplaybreaks[4]

\def \CO{{\mathcal O}}

\def \P{\mathbb P}

\def\map#1.#2.{#1 \longrightarrow #2}
\def\rmap#1.#2.{#1 \dasharrow #2}
\DeclareMathOperator{\rank}{rank}

\def\fb#1.{\underset #1 \to \times}
\def\pr#1.{\Bbb P^{#1}}
\def\ring#1.{\mathcal O_{#1}}
\def\mlist#1.#2.{{#1}_1,{#1}_2,\dots,{#1}_{#2}}

\def\uloopr#1{\ar@'{@+{[0,0]+(-4,5)} @+{[0,0]+(0,10)}
@+{[0,0]+(4,5)}}^{#1}}

\def\dloopr#1{\ar@'{@+{[0,0]+(-4,-5)} @+{[0,0]+(0,-10)}
@+{[0,0]+(4,-5)}}_{#1}}

\def\rloopd#1{\ar@'{@+{[0,0]+(5,4)} @+{[0,0]+(10,0)}
@+{[0,0]+(5,-4)}}^{#1}}

\def\lloopd#1{\ar@'{@+{[0,0]+(-5,4)} @+{[0,0]+(-10,0)}
@+{[0,0]+(-5,-4)}}_{#1}}

\long\def\ignore#1{}
\long\def\ignore#1{#1}

\begin{document}
\begin{center}
{\bf A Compactification of the Space of Holomorphic 
Maps from $\P^1$ into $\P^r$ }

\bigskip

\end{center}
\begin{center}  
{Jiayuan Lin}
\end{center}
\bigskip

{\bf {\small Abstract}} {\small Let $M_{d}(\P^r)$ be the space of $(r+1)$-tuples $(f_0,\cdots,f_r)$ modulo homothety, where $f_0,\cdots,f_r$ are homogeneous polynomials of degree $d$ in two variables. Let $M_{d}^{\circ}(\P^r)$ be the open subset of $M_{d}(\P^r)$ such that $f_0,\cdots,f_r$ have no common zeros. Then $M_{d}^{\circ}(\P^r)$ parametrizes the space of holomorphic 
maps of degree $d$ from $\P^1$ into $\P^r$. In general the boundary divisor $M_{d}(\P^r) \setminus M_{d}^{\circ}(\P^r)$ is not normal crossing. In this paper we will give a natural stratification of this boundary and show that we can process 
an iterated blow-ups along these strata (or its proper transformations) to obtain a compactification of $M_{d}^{\circ}(\P^n)$ with normal crossing divisors.} 

\section{Introduction}

Let $M_{d}(\P^r)=\P^{(d+1)(r+1)-1}$ be the space of $(r+1)$-tuples $(f_0,\cdots,f_r)$ modulo homothety, where $f_0,\cdots,f_r$ are homogeneous polynomials of degree $d$ in two variables. Let $M_{d}^{\circ}(\P^r)$ be the open subset of $M_{d}(\P^r)$ such that $f_0,\cdots,f_r$ have no common zeros. Then $M_{d}^{\circ}(\P^r)$ parametrizes the space of holomorphic maps of degree $d$ from $\P^1$ into $\P^r$. In general the boundary divisor $M_{d}(\P^r) \setminus M_{d}^{\circ}(\P^r)$ is not normal crossing, so $M_{d}(\P^r)$ is not a good compactification of the open variety $M_{d}^{\circ}(\P^r)$. In this paper we will give a natural stratification of $M_{d}(\P^r) \setminus M_{d}^{\circ}(\P^r)$ and show that an iterated blow-ups along these strata (or its proper transformations) can be carried out to obtain a compactification of $M_{d}^{\circ}(\P^n)$ with normal crossing divisors.

The existence of compactifications of an open variety by adding normal crossing divisors is guaranteed by Hironaka's Theorem on resolution of singularities. However, in practice, for certain open varieties people like to construct such compactifications explicitly. Take, for example, the famous work by W. Fulton and R. Macpherson $[1]$ on configuration spaces, and its further extensions by R. Macpherson and C. Procesi $[6]$ and Ulyanov $[7]$. Recently Y. Hu $[4]$ generalized the above results. For any given open variety meeting certain conditions, he proved that a compactification with normal crossing divisors can be obtained by blowing up along an arrangement of its subvarieties. As one of applications of his theorem, in the introduction and section $6$ of $[4]$, Y. Hu proposed the following
     
\begin{thm} A compactification of $M_{d}^{\circ} (\P^r)$ with normal 
crossing divisors can be obtained by iterated blow-ups along the strata (or its proper transformations) of $M_{d}(\P^r) \setminus M_{d}^{\circ}(\P^r)$.
\end{thm}

S. Keel observed that some problem might arise in the iterated blow-ups. The lowest stratum is smooth, along which the blow-up can be carried out. However, starting from the second one, each stratum has singularities along the lower strata. So the second blow-up can be done only if one can show that the singularities of the second stratum have been resolved after the first blow-up. Similarly in order to carry out the $k$-th $(2 \le k \le d)$ blow-up, one has to show that the singularities of the $k$-th stratum are resolved after the first $(k-1)$ blow-ups. In this 
paper we will confirm that it is the case. As a consequence, Theorem $1.1$ holds.

This paper is organized as follows. Section $2$ is a preliminary, in which we will describe a natural stratification of $M_{d}(\P^r) \setminus M_{d}^{\circ}(\P^r)$ and its related properties. 
In section $3$, we will first set up some notations, and then give a proof of Theorem $1.1$.

\medskip

{\bf Acknowledgements} I would like to thank Professor Y. Hu for his interest in and comments on this work. I am also grateful to Professor J. M$^c$Kernan and Professor V. Alexeev for helpful discussions and suggestions.

\section{preliminary}
Let $f_{i}(x,y)=\overset{d}{\underset{j=0}{\sum}} s_{i j}x^{d-j} y^{j} (i=0,1,\cdots,r)$ be $r+1$ homogeneous polynomials of degree $d$ in two variables $x$ and $y$. Let $R_k$ be the subset in $M_{d}(\P^r)$ which parametrizes $f_0,\cdots, f_r$ with
at least $d-k+1$ common roots (counting multiplicities). Then $R_1 \subset R_2 \subset \cdots \subset R_d=M_{d}(\P^r) \setminus M_{d}^{\circ}(\P^r)$ gives a natural stratification of $M_{d}(\P^r) \setminus M_{d}^{\circ}(\P^r)$.

Let $A_{(r+1)k \times (d+k)}$ be the following matrix:   

\begin{equation}
\begin{split}
A_{(r+1)k \times (d+k)} =\begin{pmatrix} 
s_{00} & s_{01} &\cdots& s_{0d}& \cdots & \cdots&\cdots&\cdots&\cdots&0\cr
s_{10} & s_{11}  &\cdots& s_{1d}&\cdots & \cdots&\cdots&\cdots&\cdots&0\cr
\cdots&\cdots&\cdots&\cdots&\cdots&\cdots&\cdots&\cdots&\cdots&\cdots\cr
s_{r0} & s_{r1}  &\cdots& s_{rd}&\cdots &\cdots& \cdots&\cdots&\cdots&0\cr
0&s_{00} & s_{01} &\cdots& s_{0d}& \cdots & \cdots&\cdots&\cdots&0\cr
0&s_{10} & s_{11}  &\cdots& s_{1d}&\cdots & \cdots&\cdots&\cdots&0\cr
\cdots&\cdots&\cdots&\cdots&\cdots&\cdots&\cdots&\cdots&\cdots&\cdots\cr
0&s_{r0} & s_{r1}  &\cdots& s_{rd}&\cdots &\cdots& \cdots&\cdots&0\cr
\cdots&\cdots&\cdots&\cdots&\cdots&\cdots&\cdots&\cdots&\cdots&\cdots\cr
\cdots & \cdots& \cdots& \cdots&\cdots&0&s_{00} & s_{01} &\cdots& s_{0d}\cr
\cdots & \cdots&\cdots&\cdots&\cdots&0&s_{10} & s_{11}  &\cdots& s_{1d} \cr
\cdots&\cdots&\cdots&\cdots&\cdots&\cdots&\cdots&\cdots&\cdots&\cdots\cr
\cdots &\cdots& \cdots&\cdots&\cdots&0&s_{r0} & s_{r1}  &\cdots& s_{rd}\cr
\end{pmatrix}_{(r+1)k\times (d+k)}
\end{split}
\end{equation}

Then the following lemmas are true.
\begin{lem} The polynomials $f_0, \cdots, f_r$ have 
at least $d-k+1$ common roots (counting multiplicities) in $\P^1$ if and only if the matrix $A_{(r+1)k \times (d+k)}$ has rank less than $2k$. 
\end{lem}
\begin{proof} See proposition $3$ in $[5]$. 
\end{proof}
\begin{lem} $R_k$ has dimension $d+kr$
\end{lem}
\begin{proof} $R_k$ parametrizes all the tuples $(f_0,\cdots,f_r)$, or equivalently, points $[s_{00}:\cdots:$

\noindent $s_{0d}:\cdots:s_{r0}:\cdots:s_{0d}]$ in $\P^{(d+1)(r+1)-1}$, such that $f_{i}(x,y)=
\overset{d}{\underset{j=0}{\sum}} s_{i j}x^{j} y^{d-j}$  $(i=0,1,\cdots,r)$ have
at least $d-k+1$ common roots (counting multiplicities). So, to compute its dimension,
it is enough to compute the dimension of the intersection of $R_k$ and the affine open sets $s_{i0} \ne 0$ $(i=0, 1,\cdots, r)$.
For each $i$, $s_{i0}^{-1} f_{i} (x,y)$ can be expressed as a product of linear factors $x+ \alpha_{ij} y (j=1,\cdots, d)$. Hence we have $(r+1)$ ways to choose $s_{i0}$ $(i=0, \cdots, r)$, $(d-k+1)$ ways to choose the common roots (counting multiplicities),  and $(k-1)$ ways to choose the rest linear factors appearing in each $s_{i0}^{-1} f_{i} (x,y)$ for $i=0,1,\cdots,r$. Therefore the dimension of the intersection of $R_k$ and the affine open sets $s_{i0} \ne 0$ $(i=0, \cdots, r)$ is $(r+1)+(d-k+1)+(k-1)(r+1)-1=d+kr$, so is that of $R_k$.   
\end{proof}

Now let us determine $R_{k}$'s smooth and singular loci. The following result is 
true.
\begin{thm} The singular locus of $R_k$ is exactly $R_{k-1}$.
\end{thm}
\begin{proof} By Lemma $2.1$, a point $[s_{00}: \cdots: s_{0d}: \cdots: s_{r0}: \cdots: s_{rd}] \in 
R_k$ if and only if the rank of the matrix  $A_{(r+1)k \times (d+k)}$ 
is less than $2k$, in other words, $R_k$ is the common zero locus of 
all the $2k \times 2k$ minor determinants of matrix $A_{(r+1)k \times (d+k)}$.

The singular locus of $R_k$ is all the points where the rank of the matrix consisting of all the derivatives of locally defining functions of $R_k$ with respect to local coordinates fails to be equal to the codimension of $R_k$ in $\P^{(d+1)(r+1)-1}$.

For any point $[s_{00}: \cdots: s_{0d}: \cdots: s_{r0}: \cdots: s_{rd}] \in 
R_{k-1}$, the rank of the matrix

\noindent $A_{(r+1)(k-1) \times (d+k-1)}$ 
is less than $2(k-1)$. By the only lemma in $[5]$, $\rank A_{(r+1)k \times (d+k)}=\rank$
\noindent $A_{(r+1)(k-1) \times (d+k-1)}+1< 2(k-1)+1=2k-1$. Suppose $s_{i_0j_0} \ne 0$ for some $i_0$ and $j_0$. Then $s_{i_0j_0}^{-2k} |M|$ are all the locally defining functions of $R_k$, where $|M|$ runs over all the $2k \times 2k$ minor determinants of matrix $A_{(r+1)k \times (d+k)}$. The derivatives of $s_{i_0j_0}^{-2k}|M|$ with respect to the local coordinate $\frac{s_{ij}}{s_{i_0j_0}}$ are linear combinations of $(2k-1) \times (2k-1)$ minor determinants of $s_{i_0j_0}^{-2k}|M|$. Since $\rank A_{(r+1)k \times (d+k)}<2k-1$, it follows that any $(2k-1) \times (2k-1)$ minor determinant of $M$ vanishes, so do the derivatives of $s_{i_0j_0}^{-2k}|M|$. Therefore all the derivatives of the local defining functions of $R_k$ vanish on $R_{k-1}$, which implies that $R_{k-1}$ is a subset of the singular locus of $R_k$.

In order to show that $R_{k-1}$ is exactly the singular locus of $R_k$, we need to show that
$R_{k}\setminus R_{k-1}$ is smooth.

Note that $R_k$ parametrizes $(r+1)$-tuples $(f_0,\cdots, f_r)$ of homogeneous polynomials of degree $d$ with at least $d-k+1$ common roots (count multiplicities). So given a point in $R_k$, we can find a homogeneous polynomial $p(x,y)$ of degree $d-k+1$ and a $(r+1)$-tuples $(g_0,\cdots,g_r)$ of homogeneous polynomials of degree $k-1$ such that $f_i(x,y)=p(x,y) g_i (x,y)$. Conversely, given any homogeneous polynomials $p(x,y)$ of degree $d-k+1$ and a $(r+1)$-tuples $(g_0,\cdots,g_r)$ of homogeneous polynomials of degree $k-1$, we can associate them with a point in $R_k$ by letting $f_i(x,y)=p(x,y) g_i( x,y)$. This correspondence allows us to define a morphism $\Phi_k: \map \P^{(r+1)k-1} \times \P^{d-k+1}. \P^{(d+1)(r+1)-1}.$ as $\Phi_k([\mu_{00}:\cdots: \mu_{0,k-1}:\mu_{10}:\cdots:\mu_{1,k-1}:\cdots:\mu_{r0}:\cdots:$
$\mu_{r,k-1}];[\nu_0: \cdots: \nu_{d-k+1}])=[ \mu_{00} \nu_0: \mu_{00} \nu_1+\mu_{01} \nu_0:\cdots: \mu_{0,k-1}\nu_{d-k+1}: \cdots: \mu_{r0} \nu_0: \mu_{r0} \nu_1+$ 
$\mu_{r1} \nu_0:$
$\cdots: \mu_{r,k-1}\nu_{d-k+1}]$, or equivalently, $\Phi_k (g_0,\cdots,g_r;p)=(g_0p,\cdots,g_r p)$, where $\P^{(r+1)k-1}$ parametrizes $(r+1)$-tuples $(g_0,\cdots,g_r)$ with $g_i (x, y)= $
$\overset{k-1}{\underset{j=0}{\sum}} \mu_{ij}$
$x^{k-1-j} y^{j}$ $(i=0,1,\cdots,r)$ and $\P^{d-k+1}$ parametrizes $p(x,y)=\overset{d-k+1}{\underset{j=0}{\sum}} \nu_{i}x^{d-k+1-i} y^{i}$. It is easy to see that $\Phi_k$ is surjective over $R_k$ and bijective (in fact isomorphic) from the preimage of $R_{k}\setminus R_{k-1}$ onto $R_{k}\setminus R_{k-1}$. Therefore $R_{k}\setminus R_{k-1}$ is smooth and $R_{k-1}$ is the singular locus of $R_k$. 
\end{proof}

We need the following lemma for future use.
\begin{lem} The rank of $A_{(r+1)k \times (d+k)}$ is less than $2k$ if and only if any its submatrix with the form
\begin{equation}
\begin{split}
B_{2k \times (d+k)}= 
\begin{pmatrix} 
* & *&\cdots& *& \cdots & \cdots&\cdots&\cdots&\cdots&0\cr
* & * &\cdots& *&\cdots & \cdots&\cdots&\cdots&\cdots&0\cr
0&* & * &\cdots& *& \cdots & \cdots&\cdots&\cdots&0\cr
0&* & * &\cdots& *&\cdots & \cdots&\cdots&\cdots&0\cr
\cdots&\cdots&\cdots&\cdots&\cdots&\cdots&\cdots&\cdots&\cdots&\cdots\cr
\cdots&\cdots&\cdots&\cdots&\cdots&\cdots&\cdots&\cdots&\cdots&\cdots\cr
\cdots & \cdots& \cdots& \cdots&\cdots&0&* & * &\cdots& *\cr
\cdots & \cdots&\cdots&\cdots&\cdots&0&* & *  &\cdots& * \cr
\end{pmatrix}_{2k \times (d+k)}
\end{split}
\end{equation}

has rank less than $2k$, where the position marked by $(*,*,\cdots,*)$ in each row is filled with one of the vectors from $(s_{i0}, s_{i1}, \cdots, s_{id}), i=0, \cdots,r$. 
\end{lem}
\begin{proof} The \lq\lq only if" part is trivial, so we only need to show the \lq\lq if" part.

We use induction on $k$.

$k=1$ is trivial.

Suppose that our lemma is true for $k-1$, i.e., that any submatrix in $A_{(r+1)(k-1) \times (d+k-1)}$ with the form $B_{2(k-1) \times (d+k-1)}$ has rank less than $2(k-1)$ implies the rank of 

\noindent $A_{(r+1)(k-1) \times (d+k-1)}$ is less than $2(k-1)$. We need to show that the rank of $A_{(r+1)k \times (d+k)}$ is less tahn $2k$ if all it submatrices with the form $B_{2k \times (d+k)}$ 
have rank less than $2k$. 

Two subcases:

Subcase I: All of the submatrices with the form $B_{2(k-1) \times (d+k-1)}$ have rank less than $2(k-1)$. By our inductive assumption 
$ \rank A_{(r+1)(k-1) \times (d+k-1)}<2(k-1)$. By the only lemma in $[5]$, which says that $\rank A_{(r+1)k \times (d+k)}= \rank A_{(r+1)(k-1) \times (d+k-1)}+1$ if $\rank A_{(r+1)(k-1) \times (d+k-1)}<2(k-1)$, we have that $\rank A_{(r+1)k \times (d+k)}= $

\noindent $\rank A_{(r+1)(k-1) \times (d+k-1)}+1< 2(k-1)+1<2k$.

Subcase II: At least one of the submatrices with the form $B_{2(k-1) \times (d+k-1)}$ has rank equal to $2(k-1)$. We will show that $\rank A_{(r+1)k \times (d+k)}=2k-1<2k$.

Denote the submatrix consisting of vectors from the $[(r+1)\cdot(i-1)+1]$-th row to the $[(r+1)\cdot i]$-th row in $A_{(r+1)k \times (d+k)}$ the $i$-th block, where $i=1, \cdots, k$. By the assumption of subcase II, there is at least one matrix with the form $B_{2(k-1) \times (d+k-1)}$ and rank $2(k-1)$. Pick a such matrix and still denote it $B_{2(k-1) \times (d+k-1)}$.   
Let $\varepsilon_i, \eta_i, i=1, \cdots,k-1$ be the pairs of vectors obtained from $B_{2(k-1) \times (d+k-1)}$ by extending each row one more position from behind and filling it with zero. These vectors appear in the first $(k-1)$ blocks in $A_{(r+1)k \times (d+k)}$. There exists a vector in the $k$-th block of $A_{(r+1)k \times (d+k)}$, say $\varepsilon_k$, which is linearly independent of any set of linearly independent vectors appearing in the first $(k-1)$ blocks. In particular, $\varepsilon_k$ is linearly independent of $\varepsilon_i, \eta_i (i=1, \cdots,k-1)$. Pick a vector $\xi$ other than $\varepsilon_k$ in the $k$-th block. Then $\xi$, together with $\varepsilon_k$ and $\varepsilon_i, \eta_i (i=1, \cdots,k-1)$, forms a matrix with the form $B_{2k \times (d+k)}$. By our assumption, its rank is less than $2k$. So $\xi$ is a linear combination of $\varepsilon_k$ and $\varepsilon_i, \eta_i, i=1, \cdots,k-1$.

Let $V$ be the vector space spanned by $\varepsilon_k$ and $\varepsilon_i, \eta_i (i=1, \cdots,k-1)$. We need to show that $V$ contains any row in $A_{(r+1)k \times (d+k)}$. 

Since $\xi$ is an arbitrary vector in the $k$-th block and we have showed that it is a linear combination of $\varepsilon_k$ and $\varepsilon_i, \eta_i (i=1, \cdots,k-1)$, it follows that $V$ contains the $k$-th block.

We need the following claim before we can move on.

\begin{claim} If the rank of a submatrix $M$ of $A_{(r+1)k \times (d+k)}$ consisting of the $k$-th block and another $s$ blocks from the first $(k-1)$ ones is $2s+1$, then $V$ contains all these $s+1$ blocks. Moreover, $V$ contains all the $k$ blocks in $A_{(r+1)k \times (d+k)}$. 
\end{claim}
\begin{proof}      The pairs of $\varepsilon_i, \eta_i$ in the $s$ blocks and $\varepsilon_k$ in the $k$-th block in $A_{(r+1)k \times (d+k)}$ form a submatrix of $M$ and it has the same rank as $M$, thus any row in $M$ is a linear combination of those $\varepsilon_i, \eta_i$ ($i=1,\cdots,k-1$) and $\varepsilon_k$. Since $\varepsilon_i, \eta_i$ ($i=1,\cdots,k-1$) and $\varepsilon_k$ are vectors in $V$, it follows that $V$ contains all these $s+1$ blocks.

If $s=k-1$, then $V$ contains all the $s+1=k$ blocks in $A_{(r+1)k \times (d+k)}$, we are done.

Suppose $s<k-1$, we will show that we can produce $s'$ blocks $(s'> s)$ among the first $(k-1)$ ones such that the rank of the matrix consisting of these $s'$ blocks and the $k$-th block is $2s'+1$. The same argument as the first paragraph gives that $V$ contains all these $s'+1$ blocks. If $s'<k-1$, replace $s'$ by $s$ and repeat the above process. After finite many steps we will end at $s'=k-1$, which exactly means that the rank of $A_{(r+1)k \times (d+k)}$ is $2(k-1)+1=2k-1$ and $V$ contains all the $k$ blocks in $A_{(r+1)k \times (d+k)}$. 

Assume $V$ contains the $i_j$-th $(j=1, \cdots, s)$ blocks and the $k$-th block. Let $V_1$ be the space spanned by the rows in the $i_j$-th blocks $(j=1, \cdots, s)$. By our assumption, $\dim V_1=2s$. In fact $V_1$ is generated by $\varepsilon_{i_j}, \eta_{i_j} (j=1, \cdots,s)$. We have the following four possible cases:

(1) If $i_1<\cdots<i_s$ are consective natural numbers and $i_{s} <k-1$, then we can shift the $i_j$-th block left by one unit for all $j=1,\cdots, s$ , where \lq\lq shift left by one unit" means that the piece of $s_{i0}, s_{i1}, \cdots, s_{id}$ $(i=0, \cdots, r)$ in each row in the $i_j$-th block is shifted left by one unit. Let $V_2$ be the space spanned by the row vectors in the $(i_j+1)$-th blocks, $j=1, \cdots, s$. Then $\dim V_2= \dim V_1=2s$. There are $(s-1)$ blocks contained in both $V_1$ and $V_2$ , each with a pair of $\varepsilon_i, \eta_i$, so $\dim (V_1 \cap V_2) \ge 2(s-1)$. Therefore $\dim (V_1+V_2)= \dim V_1 +\dim V_2 -\dim (V_1 \cap V_2) \le 2s+ 2s-2(s-1)=2(s+1)$. However, $V_1+V_2$ contains $s+1$ blocks among the first $(k-1)$ ones in $A_{(r+1)k \times (d+k)}$, which is at least of rank $2(s+1)$, so $\dim (V_1+V_2)= 2(s+1)$. The vectors $\varepsilon_{i_j}, \eta_{i_j}, j=1, \cdots,s$ and $\varepsilon_{i_{s}+1}, \eta_{i_{s}+1}$ in $V_1+V_2$ are linearly independent, so it spans $V_1+V_2$. Since $\varepsilon_{i_j}, \eta_{i_j}, j=1, \cdots,s$ and $\varepsilon_{i_{s}+1}, \eta_{i_{s}+1}$ are vectors in $V$, we conclude that $V_1+V_2$ is a subspace of $V$. So the rank of the submatrix consisting of these $s+1$ blocks and the $k$-th block is $2(s+1)+1$, and $V$ contains all these $s+2$ blocks.

(2) If $i_1<\cdots<i_s$ are consective natural numbers and $i_{s} =k-1$, then we can shift the $i_j$-th block right by one unit for $j=1,\cdots, s$ and repeat the above process to get that the rank of the submatrix consisting of $s+1$ blocks among the first $(k-1)$ ones in $A_{(r+1)k \times (d+k)}$ and the $k$-th block is $2(s+1)+1$, and $V$ contains all these $s+2$ blocks.

(3) If $i_1<\cdots<i_s$ are not consective natural numbers and $i_{s} <k-1$, then we can shift the $i_j$-th block left by one unit for $j=1,\cdots, s$. Define $V_1$ and $V_2$ as before. Suppose $V_1$ and $V_2$ share $t$ blocks. Then that $i_j$ $(j=1, \cdots,s)$ are not consective natural numbers implies $t \le s-1$. Repeating the same process as in case (1), we can show that $\dim (V_1+V_2)= 2(s+s-t)=2s+2(s-t) \ge 2(s+1)$ and $V_1+V_2$ is a subspace of $V$. Moreover, the rank of the submatrix consisting of the $k$-th block and the $(2s-t)$ blocks from the first $(k-1)$ ones in $A_{(r+1)k \times (d+k)}$ is $2(2s-t)+1\ge 2(s+1)+1$, and $V$ contains all these $2s-t+1 \ge s+2$ blocks.

(4) If $i_1<\cdots<i_s$ are not consective natural numbers and $i_{s} =k-1$, we still shift the $i_j$-th block left by one unit for $j=1,\cdots, s$. Define $V_1$ and $V_2$ as before. 
Suppose $V_1$ and $V_2$ share $t$ bolcks. Since $i_j$, $j=1, \cdots,s$ are not consective natural numbers and the $i_{s}+1=k$-th block is not in $V_1$, we have $t \le s-2$. Let $L\varepsilon_{i_j}, L\eta_{i_j} (j=1, \cdots,s)$ be the vectors in $V_2$ obtained by shifting $\varepsilon_{i_j}, \eta_{i_j} (j=1, \cdots,s)$ left by one unit. Then $L\varepsilon_{i_j}, L\eta_{i_j} (j=1, \cdots,s)$ are linearly independent. The vectors $L\varepsilon_{i_{s}}, L\eta_{i_{s}}$ are in the $k$-th block. We know that any vector in the $k$-th block can be expressed as linear combination of $\varepsilon_{i_j}, \eta_{i_j} (j=1, \cdots,s)$ and $\varepsilon_k$.  Assume $L\varepsilon_{i_{s}}=a \varepsilon_k+ \overset{s}{\underset{j=1}{\sum}} l_j \varepsilon_{i_j} + \overset{s}{\underset{j=1}{\sum}} m_j \eta_{i_j}$ and $L\eta_{i_{s}}=
b \varepsilon_k+ \overset{s}{\underset{j=1}{\sum}} l'_j \varepsilon_{i_j} + \overset{s}{\underset{j=1}{\sum}} m'_j \eta_{i_j}$. If one of $a$ and $b$ is zero, then $L\varepsilon_{i_{s}}$ or $L\eta_{i_{s}}$ is in $V_1$. Since $L\varepsilon_{i_{s}}$ and $L\eta_{i_{s}}$ are also in $V_2$, it follows that at least one of $L\varepsilon_{i_{s}}$ and $L\eta_{i_{s}}$ is in $V_1 \cap V_2$. Now the $2t$ pairs of $L\varepsilon_{i_j}, L\eta_{i_j}$ in the $t$ common bolcks appear in $V_1 \cap V_2$, together with one of  $L\varepsilon_{i_{s}}$ and $L\eta_{i_{s}}$, imply that $\dim (V_1 \cap V_2) \ge 2t+1$. If neither of $a$ and $b$ is zero, then $b L\varepsilon_{i_{s}}- a L\eta_{i_{s}} \ne 0$  is in $V_1 \cap V_2$. $b L\varepsilon_{i_{s}}- a L\eta_{i_{s}}$, together with the $2t$ pairs of $L\varepsilon_{i_j}, L\eta_{i_j}$ in the $t$ common bolcks in $V_1 \cap V_2$, also gives that $\dim (V_1 \cap V_2) \ge 2t+1$. So $\dim (V_1+V_2)= \dim V_1 + \dim V_2- \dim (V_1 \cap V_2) \le 2s+2s- (2t+1)= 2s+2(s-1-t)+1$. But $V_1+V_2$ contains $s+s-t-1$ blocks from the first $(k-1)$ ones and the $k$-th block. The vectors consisting of the pairs of $\varepsilon_{i_{j}}, \eta_{i_{j}}$ and $\varepsilon_k$ generate a subspace $W$ of $V_1+V_2$ of dimension $2(2s-1-t)+1=\dim (V_1+V_2)$. Hence $W=V_1+V_2$. Since all the pairs of $\varepsilon_{i_{j}}, \eta_{i_{j}}$ and $\varepsilon_k$ are vectors in $V$, $V$ contains $V_1+V_2$. Therefore the rank of the submatrix consisting of these $s+s-t-1=s+(s-1-t)$ blocks from the first $(k-1)$ ones in $A_{(r+1)k \times (d+k)}$ and the $k$-th block is $2[s+(s-1-t)]+1 \ge 2(s+1)+1$ and $V$ contains all these $s+(s-1-t)+1 \ge s+2$ blocks. 

\end{proof}

By Claim $2.5$ we only need to show that there exists a submatrix $M$ consisting of the $k$-th block and another $s$ blocks from the first $(k-1)$ ones in 
$A_{(r+1)k \times (d+k)}$ such that $\rank M=2s+1$.

Let us assign each block in $A_{(r+1)k \times (d+k)}$ a level number.

The level $0$ consists of only one block, the $k$-th one. We have known that $V$ contains the $k$-th block.

A block among the first $(k-1)$ blocks in $A_{(r+1)k \times (d+k)}$, say the $p$-th one, is in the level $1$ if there exists a row vector $\xi$ in the $k$-th block such that $\xi= \overset{k-1}{\underset{i=1}{\sum}} a_i \varepsilon_{i} + \overset{k-1}{\underset{i=1}{\sum}} b_i \eta_{i} + a \varepsilon_k$ and at least one of $a_p$ and $b_p$ is non-zero. WLOG, say $a_p \ne 0$. Then $\varepsilon_{p}$ is a linear combination of $\varepsilon_{1},\eta_{1}, \cdots,\varepsilon_{p-1},\eta_{p-1},\eta_p,\varepsilon_{p+1},\eta_{p+1}, \cdots, \varepsilon_{k-1},\eta_{k-1}, \varepsilon_{k},\xi$. For any vector $\zeta$ in the $p$-th block, $\varepsilon_{1},\eta_{1}, \cdots,\varepsilon_{p-1},\eta_{p-1},\zeta,\eta_p,\varepsilon_{p+1},\eta_{p+1}, \cdots, \varepsilon_{k-1},\eta_{k-1}, \varepsilon_{k},\xi$ form a matrix with the form $B_{2k \times (d+k)}$. By our assumption, any such a matrix has rank less than $2k$.  So $\zeta$ is a linear combination of $\varepsilon_{1},\eta_{1}, \cdots,\varepsilon_{p-1},\eta_{p-1},\eta_p,\varepsilon_{p+1},\eta_{p+1}, \cdots, \varepsilon_{k-1},\eta_{k-1}, \varepsilon_{k},\xi$. All the vectors $\varepsilon_{1},\eta_{1}, \cdots,\varepsilon_{p-1},\eta_{p-1},\eta_p,\varepsilon_{p+1},\eta_{p+1}, \cdots, \varepsilon_{k-1},\eta_{k-1}, \varepsilon_{k},\xi$ are in $V$, so is $\zeta$. Since $\zeta$ is arbitrarily picked in the $p$-th block, it follows that 
$V$ contains the $p$-th block.

Suppose the number of all blocks in the level $1$ is $s_1$, then $s_1 \ge 1$. Since $V$ contains all the blocks in the levels $0$ and $1$, we are done if $s_1=k-1$. Otherwise $s_1<k-1$. From the choice of the level $1$, it is easy to see that the vectors $\varepsilon_{i_1},\eta_{i_1}, \cdots, \varepsilon_{i_{s_1}},\eta_{i_{s_1}},\varepsilon_k$ span the $k$-th block, where $i_j (j=1,\cdots,s_1)$ are the indices of the $s_1$ blocks in the level $1$. If these vectors also span the $s_1$ blocks, we get a submatrix satisfying the assumption in Claim $2.5$, hence we are done by Claim $2.5$. Otherwise there exists at least one row $\zeta$ in a block in the level $1$ , say the $p$-th one, such that the coefficients of $\varepsilon_{q}$ or $\eta_{q}$ is non-zero for some $q$ if we express $\zeta$ as a linear combination of $ \varepsilon_{1},\eta_{1}, \cdots, \varepsilon_{k-1},\eta_{k-1},\varepsilon_k$, where $ 1\le q \le k-1$ and the $q$-th block is not in the levels $0$ and $1$.  We will show that $V$ contains this $q$-th block.

To show that $V$ contains the $q$-th block, we apply the same trick as we did before in proving that $V$ contains the $p$-th block in the level $1$.

Since the $p$-th block is in the level $1$, there exists $\xi$ in the $k$-th block such that $\xi= \overset{k-1}{\underset{i=1}{\sum}} a_i \varepsilon_{i} + \overset{k-1}{\underset{i=1}{\sum}} b_i \eta_{i} + a \varepsilon_k$ and at least one of $a_p$ and $b_p$ is non-zero. WLOG, we assume $a_p \ne 0$. From the choice of the level $1$, we see that $\xi$ is in fact a linear combination of the pairs of $\varepsilon_{i_1},\eta_{i_1},\cdots,\varepsilon_{i_{s_1}},\eta_{i_{s_1}}$ and $\varepsilon_k$, where $i_j$ $(j=1, \cdots, s_1)$ are the indices of the blocks in the level $1$. Suppose $i_x=p$ for some $1 \le x \le s_1$. It is easy to see that the vectors $\{\varepsilon_{i_1},\eta_{i_1},\cdots, \varepsilon_{i_{s_1}},\eta_{i_{s_1}},\varepsilon_k\}$ are linearly equivalent to $\{\varepsilon_{i_1},\eta_{i_1},\cdots,\varepsilon_{i_{x-1}},\eta_{i_{x-1}}$
$,\eta_{i_x}=\eta_p,\varepsilon_{i_{x+1}},\eta_{i_{x+1}}, \cdots, \varepsilon_k, \xi\}$. 
Obviously, the vectors $\{\varepsilon_{1},\eta_{1},\cdots,\varepsilon_{p-1},\eta_{p-1}, \eta_p, \varepsilon_{p+1},\eta_{p+1},\cdots,$  

\noindent $\varepsilon_{k-1},\eta_{k-1},\varepsilon_k,\xi\}$ are linearly equivalent to $\{\varepsilon_{1},\eta_{1},\cdots,\cdots,\varepsilon_{k-1},\eta_{k-1},\varepsilon_k,\}$, so they form a new basis of $V$. The equivalence between the vectors $\{\varepsilon_{i_1},\eta_{i_1},\cdots, \varepsilon_{i_{s_1}},\eta_{i_{s_1}},\varepsilon_k\}$ and $\{\varepsilon_{i_1},\eta_{i_1},\cdots,\varepsilon_{i_{x-1}},\eta_{i_{x-1}}$
$,\eta_{i_x}=\eta_p,\varepsilon_{i_{x+1}},\eta_{i_{x+1}}, \cdots, \varepsilon_k, \xi\}$ implies that the coefficients of $\varepsilon_{q}$ or $\eta_{q}$ is also non-zero when we expand $\zeta$ under this new basis. WLOG, assume the coefficients of $\varepsilon_{q}$ is not zero. Then clearly
$\{\varepsilon_{1},\eta_{1},\cdots,\varepsilon_{p-1},\eta_{p-1}, \zeta, \eta_p, \varepsilon_{p+1},\eta_{p+1},\cdots,\varepsilon_{q-1},\eta_{q-1},$
$\eta_q, \varepsilon_{q+1},\eta_{q+1},\cdots,\varepsilon_{k-1},$
$\eta_{k-1},\varepsilon_k,\xi\}$ form a basis of $V$. Any vector in the $q$-th block, together with this basis, forms a matrix with the form  $B_{2k \times (d+k)}$. By our assumption, the rank of such a matrix is less than $2k$. So any vector in the $q$-th block is a linear combination of $\{\varepsilon_{1},\eta_{1},\cdots,\varepsilon_{p-1},\eta_{p-1}, \zeta, \eta_p, \varepsilon_{p+1},\eta_{p+1},\cdots,\varepsilon_{q-1},\eta_{q-1}, \eta_q, \varepsilon_{q+1},\eta_{q+1},\cdots,\varepsilon_{k-1},\eta_{k-1},\varepsilon_k,$

\noindent $\xi\}$. Therefore $V$ contains the $q$-th block.

Let the level $2$ consisting of all the possible $q$-th block like the above one. Denote $s_2$ the total number of the blocks in the level $2$. Then $V$ contains all the blocks in the levels $0$, $1$ and $2$. If one of the following two cases happens, we are done.

(a) All the blocks in $A_{(r+1)k \times (d+k)}$ have already appeared in the levels $0$, $1$ or $2$, then $\rank A_{(r+1)k \times (d+k)}=\dim V=2k-1$. 

or 

(b) The vectors $\varepsilon_{i_1},\eta_{i_1},\cdots,\varepsilon_{i_{s_1}},\eta_{i_{s_1}},\varepsilon_{i_{s_1+1}},
\eta_{i_{s_1+1}},\cdots,
\varepsilon_{i_{s_1+s_2}},\eta_{i_{s_1+s_2}}$ and $\varepsilon_k$ span the level $2$, where $i_j, j=1, \cdots,s_1$ are the indices of the blocks in the level $1$ and $i_j, j=s_1+1, \cdots,s_1+s_2$ are the indices of the blocks in the level $2$ . 

If the case $(b)$ happens, then the blocks in the levels $0, 1$ and $2$ form a submatrix of rank $2(s_1+s_2)+1$. By Claim $2.5$, $V$ contains all the blocks in $A_{(r+1)k \times (d+k)}$, so $\rank A_{(r+1)k \times (d+k)}$

\noindent $=\dim V=2k-1$.

If neither case above happens, then we can proceed to pick up blocks in the level $3$. 

Since there are only $k$ blocks in total in $A_{(r+1)k \times (d+k)}$, after finitely many steps, say $r$, we must either end at the case that $V$ contains all the blocks in the levels $0, \cdots, r$ and these levels have already included all the blocks in $A_{(r+1)k \times (d+k)}$, or the vectors $\varepsilon_{i_1},\eta_{i_1},\cdots,\varepsilon_{i_{s_1}},\eta_{i_{s_1}},$
$\varepsilon_{i_{s_1+1}},
\eta_{i_{s_1+1}},\cdots,
\varepsilon_{i_{s_1+s_2}},\eta_{i_{s_1+s_2}},\cdots, \varepsilon_{i_{s_1+\cdots+s_r}},\eta_{i_{s_1+\cdots+s_r}}$ and $\varepsilon_k$ span the level $r$, where $i_j (j$
$=s_1+\cdots+s_{p-1}+1, \cdots,s_1+\cdots+s_{p-1}+s_p)$ are the indices of the blocks in the level $p$, ($p=1, \cdots,r$). 
In either case we have that the rank of $A_{(r+1)k \times (d+k)}=\dim V=2k-1<2k$.  

Sharp eyes may have detected that there is a missing point in the above argument, that is, we did not show that $V$ contains all the blocks in the level $l (1 \le l \le r)$. This can be done by induction on $l$. 

The following property $(P_l)$ is crucial and it is true for any level $l (1 \le l \le r)$.

$(P_l)$: For any block, say $\bar{q}$-th, in the level $l$, there exists a set of linearly independent vectors $\{\varepsilon'_{k}, \eta'_{k}, \varepsilon'_{i_1},\eta'_{i_1},\cdots,\varepsilon'_{i_{s_1}},\eta'_{i_{s_1}},\varepsilon'_{i_{s_1+1}},
\eta'_{i_{s_1+1}},\cdots,
\varepsilon'_{i_{s_1+s_2}},\eta'_{i_{s_1+s_2}},\cdots, \varepsilon'_{i_{s_1+\cdots+s_{l-1}+1}},$

\noindent $\eta'_{i_{s_1+\cdots+s_{l-1}+1}},\cdots,\varepsilon'_{i_{s_1+\cdots+s_l}},\eta'_{i_{s_1+\cdots+s_l}}\} \setminus \{\varepsilon'_{\bar{q}}\}$, which is linearly equivalent to the set of vectors $\{\varepsilon_{i_1},\eta_{i_1},\cdots,\varepsilon_{i_{s_1}},\eta_{i_{s_1}},\varepsilon_{i_{s_1+1}},
\eta_{i_{s_1+1}},\cdots,
\varepsilon_{i_{s_1+s_2}},\eta_{i_{s_1+s_2}},\cdots, \varepsilon_{i_{s_1+\cdots+s_{l-1}+1}},\eta_{i_{s_1+\cdots+s_{l-1}+1}},\cdots,$

\noindent $\varepsilon_{i_{s_1+\cdots+s_l}},\eta_{i_{s_1+\cdots+s_l}},\varepsilon_k\}$, where  
$i_j (j=s_1+\cdots+s_{p-1}+1, \cdots,s_1+\cdots+s_{p-1}+s_p)$ are the indices of the blocks in the level $p$ ($p=1, \cdots,l$) and $\varepsilon'_{k}, \eta'_{k}$ are vectors in the $k$-th block (the level $0$).

$(P_1)$ has been shown true in the above.

Suppose $(P_{l-1})$ is true.

For any block, say $\bar{q}$-th, in the level $l$, there exists at least one row $\zeta$ in a block in the level $l-1$, say the $\bar{p}$-th one, such that at least one of the coefficients of $\varepsilon_{\bar{q}}$ and $\eta_{\bar{q}}$ is non-zero when we express $\zeta$ as a linear combination of $\{ \varepsilon_{1},\eta_{1}, \cdots, \varepsilon_{k-1},\eta_{k-1}, \varepsilon_k\}$ and no row in the levels less than $l-1$ has this property. Since $(P_{l-1})$ is true, that is, there exists a set of linearly independent vectors $\{\varepsilon'_{k}, \eta'_{k}, \varepsilon'_{i_1},\eta'_{i_1},\cdots,\varepsilon'_{i_{s_1}},\eta'_{i_{s_1}},$
$\varepsilon'_{i_{s_1+1}},
\eta'_{i_{s_1+1}},\cdots,
\varepsilon'_{i_{s_1+s_2}},\eta'_{i_{s_1+s_2}},\cdots,$

\noindent $\varepsilon'_{i_{s_1+\cdots+s_{l-1}}},\eta'_{i_{s_1+\cdots+s_{l-1}}}\} \setminus \{\varepsilon'_{\bar{p}}\}$, which is linearly equivalent to the set of vectors $\{\varepsilon_{i_1},\eta_{i_1},\cdots,$

\noindent $\varepsilon_{i_{s_1}},\eta_{i_{s_1}},$
$\varepsilon_{i_{s_1+1}},
\eta_{i_{s_1+1}},\cdots,
\varepsilon_{i_{s_1+s_2}},\eta_{i_{s_1+s_2}},\cdots, \varepsilon_{i_{s_1+\cdots+s_{l-1}}},\eta_{i_{s_1+\cdots+s_{l-1}}},\varepsilon_k\}$, where  
$i_j (j=s_1+\cdots+s_{p-1}+1, \cdots,s_1+\cdots+s_{p-1}+s_p)$ are the indices of the blocks in the level $p$ ($p=1, \cdots,l-1$) and $\varepsilon'_{k}, \eta'_{k}$ are vectors in the $k$-th block. 
The linear equivalence between these two sets of vectors implies that $\zeta$  still has at least one non-zero coefficient in front of $\varepsilon_{\bar{q}}$ or $\eta_{\bar{q}}$ when we express it as a linear combination of $\{\varepsilon'_{k}, \eta'_{k}, \varepsilon'_{i_1},\eta'_{i_1},\cdots,\varepsilon'_{i_{s_1}},\eta'_{i_{s_1}},$
$\varepsilon'_{i_{s_1+1}},
\eta'_{i_{s_1+1}},\cdots,
\varepsilon'_{i_{s_1+s_2}},\eta'_{i_{s_1+s_2}},$
$\cdots,\varepsilon_{i_{s_1+\cdots+s_{l-1}+1}},\eta_{i_{s_1+\cdots+s_{l-1}+1}},\cdots,$
$\varepsilon_{i_{s_1+\cdots+s_l}},\eta_{i_{s_1+\cdots+s_l}}\} $

\noindent $\setminus \{\varepsilon'_{\bar{p}}\}$. WLOG, assume the coefficient of $\varepsilon_{\bar{q}}$ is non-zero. Then  $\{\varepsilon'_{k}, \eta'_{k}, \varepsilon'_{i_1},\eta'_{i_1},\cdots,$
$\varepsilon'_{i_{s_1}},\eta'_{i_{s_1}},$
$\varepsilon'_{i_{s_1+1}},
\eta'_{i_{s_1+1}},\cdots,
\varepsilon'_{i_{s_1+s_2}},\eta'_{i_{s_1+s_2}},$
$\cdots$
$,\varepsilon_{i_{s_1+\cdots+s_{l-1}+1}},\eta_{i_{s_1+\cdots+s_{l-1}+1}},\cdots,$
$\varepsilon_{i_{s_1+\cdots+s_l}},\eta_{i_{s_1+\cdots+s_l}}, \zeta\} \setminus \{\varepsilon'_{\bar{p}}, \varepsilon'_{\bar{q}}\}$ is a set of linearly independent vectors, which is with the form required in $(P_l)$ and linearly equivalent to the set of vectors $\{\varepsilon_{i_1},\eta_{i_1},\cdots,\varepsilon_{i_{s_1}},\eta_{i_{s_1}},\varepsilon_{i_{s_1+1}},
\eta_{i_{s_1+1}},\cdots,
\varepsilon_{i_{s_1+s_2}},$

\noindent $\eta_{i_{s_1+s_2}} ,\cdots, \varepsilon_{i_{s_1+\cdots+s_{l-1}+1}}, \eta_{i_{s_1+\cdots+s_{l-1}+1}}
,\cdots,$
$\varepsilon_{i_{s_1+\cdots+s_l}},\eta_{i_{s_1+\cdots+s_l}},\varepsilon_k\}$. Hence  
$(P_l)$ is true.

It is easy to see that $\{\varepsilon'_{k}, \eta'_{k}, \varepsilon'_{i_1},\eta'_{i_1},\cdots,\varepsilon'_{i_{s_1}},\eta'_{i_{s_1}},$
$\varepsilon'_{i_{s_1+1}},
\eta'_{i_{s_1+1}},\cdots,
\varepsilon'_{i_{s_1+s_2}},\eta'_{i_{s_1+s_2}},$
$\cdots$
$,\varepsilon_{i_{s_1+\cdots+s_{l-1}+1}},\eta_{i_{s_1+\cdots+s_{l-1}+1}},\cdots,$
$\varepsilon_{i_{s_1+\cdots+s_l}},\eta_{i_{s_1+\cdots+s_l}}, \zeta\} \setminus \{\varepsilon'_{\bar{p}}, \varepsilon'_{\bar{q}}\}$, together with all the other $\varepsilon_i,\eta_i$ not in the levels $0, 1, \cdots,$ and $l$, form a basis of $V$. Any row in the $\bar{q}$-th block, together with this basis, form a matrix with the form $B_{2k \times (d+k)}$. So this row is a linear combination of the above basis, hence it is in $V$. Since we pick the row arbitrarily from any $\bar{q}$-th block in the level $l$, it follows that $V$ contains the level $l$. This completes
the proof of the subcase II. Hence the \lq\lq if" part in the lemma $2.4$ is true. As we said before, the \lq\lq only if" part is trivial, so lemma $2.4$ follows.

\end{proof}

\section{Iterated Blow-ups}

In this section we will prove the theorem $1.1$.

We follow J. Harris $[2]$ on the definition of blow-up.

\begin{defn} 
Let $X\subseteq \P^{m}$ be any projective variety and $Y\subseteq X$ a subvariety. 
Take a collection $F_0, \cdots, F_n$ of homogeneous polynomials in $I(X)$ of the 
same degree generating an ideal with saturation $I(Y)$ (this does not have to be 
a minimal set). Consider the rational map 
$$\varphi: \rmap X.\P^{n}.$$
\noindent given by 
$$\varphi(x)=[F_0, \cdots, F_n]$$
\noindent Clearly, $\varphi$ is regular on the complement $X \setminus Y$, and in
general won't be on $Y$; thus the graph $\Gamma_{\varphi}$ (closure of the graph 
of the rational map $\varphi$) will map isomorphically to $X$ away from $Y$, but 
not in general over $Y$. The graph $\Gamma_{\varphi}$, together with the 
projection $\pi: \map \Gamma_{\varphi}. X.$, is called the blow-up of $X$ along 
$Y$ and sometimes denoted as $Bl_{Y}X$ or simply $\tilde{X}$. The inverse image 
$E=\pi^{-1}(Y) \subseteq \tilde{X}$ is called the exceptional divisor. 
\end{defn}

Note that this definition of blow-up is a generalization of the usual one in Theorem $1.1$, which requires both $X$ and $Y$ are smooth. In the case of $X$ and $Y$ both smooth, the above definition coincides with the usual one. 

Since we start the iterated blow-ups from a smooth variety $\P^{(d+1)(r+1)-1}$, to carry out the iterated process, we only need to show that the center we blow up along at each step is smooth. The advantage of using this equivalent definition is that we can construct blow-ups  explicitly, while the disadvantage is that it brings us heavy notations.

We will prove our theorem by induction on $d$.

$d=1$.

There is only one stratum $R_1$, which is smooth and of dimension $1+r$. Blow up $\P^{2(r+1)-1}$ along $R_1$ produces a smooth variety. Its boundary  consists of only one divisor, the exceptional one, which is of course normal crossing. Hence Theorem $1.1$ is trivially true for $d=1$.

Suppose that Theorem $1.1$ is true for all the spaces of holomorphic maps from $\P^1$ into $\P^r$ of degree less than $d$. We need to prove that it is also true for the space of holomorphic maps from $\P^1$ into $\P^r$ of degree $d$.

For the space of holomorphic maps from $\P^1$ into $\P^r$ of degree $d$, its boundary $M_{d}(\P^r) \setminus M_{d}^{\circ}(\P^r)$ in $\P^{(d+1)(r+1)-1}$ has a natural stratification $R_1 \subset R_2 \cdots \subset R_k$. $R_1$ is smooth, so we can blow up $\P^{(d+1)(r+1)-1}$ along $R_1$ to get $\Gamma_{\varphi_1}$, where $\varphi_1$ is a rational map associated with the ideal $I(R_1)$. However, $R_2$ is singular along $R_1$, in order to carry out the second blow-up in the usual sense, we need to show that the proper transformation $\widetilde{R}_2$ of $R_2$ in $\Gamma_{\varphi_1}$ is smooth. Similarly, assume that we have carried out the first $k$ blow-ups. If $k=d$, we are done. If $k<d$, to proceed to the next blow-up, we have to show that the proper transformation $\widetilde{R}_{k+1}$ of $R_{k+1}$ in $\Gamma_{\varphi_k}$ is smooth.

The idea to show that $\widetilde{R}_{k+1}$ is smooth is very simple. We first construct a birational morphism $\Phi_{k+1}: \map \P^{(k+1)(r+1)-1} \times \P^{d-k}. R_{k+1} \subset \P^{(d+1)(r+1)-1}.$. The space $\P^{(k+1)(r+1)-1}$ can be viewed as $M_{k}(\P^r)$. Since $k<d$, by our inductive assumption, we can carry out the iterated blow-ups to get a compactification of $M_{k}^{\circ}(\P^r)$. Denote the final outcome $\Gamma'_{\varphi_k}$, which is a smooth variety. We show that there is an isomorphism between  $\Gamma'_{\varphi_k} \times \P^{d-k}$ and $\widetilde{R}_{k+1}$. Therefore $\widetilde{R}_{k+1}$ is smooth.

\begin{rem} The reader may feel our induction a little weird. It seems that we use little information on the first $k$ blow-ups in the proof of the smoothness of $\widetilde{R}_{k+1}$. However, we do need to assume that we can carry them out to get a smooth variety $\Gamma_{\varphi_k}$ before we can move ahead.  
\end{rem}
  
Now let us show that the proper transformation $\widetilde{R}_{2}$ of $R_2$ in $\Gamma_{\varphi_1}$ is smooth.

Consider the birational morphism $\Phi_2: \map \P^{2(r+1)-1} \times \P^{d-1}. R_2 \subset \P^{(d+1)(r+1)-1}.$, where $\Phi_2([\mu_{00}:\mu_{01}:\mu_{10}:\mu_{11}:$
$\cdots: \mu_{r0}:\mu_{r1}]; [\nu_0: \cdots: \nu_{d-1}])=[ \mu_{00} \nu_0: \mu_{00} \nu_1+\mu_{01} \nu_0: \cdots:\mu_{00} \nu_{d-1}+\mu_{01} \nu_{d-2}: \mu_{01}\nu_{d-1}: \cdots: \mu_{r0} \nu_0: \mu_{r0} \nu_1+$ 
$\mu_{r1} \nu_0:\cdots:\mu_{r0} \nu_{d-1}+\mu_{r1} \nu_{d-2}:$

\noindent $\mu_{r1}\nu_{d-1}]$. Since $\P^{2(r+1)-1}$ parametrizes the space of $(r+1)$-tuples polynomials of degree $1$ modulo homothety, we can blow it up along its unique stratum. Denote the associated rational map $\varphi'_1$ and the resulting variety $\Gamma_{\varphi'_1}$. Then we claim that there exists a morphism $F$ from $\Gamma_{\varphi'_1} \times \P^{d-1}$ into $\widetilde{R}_2 \subset \P^{(d+1)(r+1)-1}$, which is an isomorphism. So $\widetilde{R}_2$ is smooth.

By lemma $2.4$ the unique stratum $R'_1$ in $\P^{2(r+1)-1}$ is generated by all $2 \times 2$ determinants $|\begin{smallmatrix}\mu_{i0} &\mu_{i1}\cr \mu_{j0}&\mu_{j1}\cr
\end{smallmatrix}|$, where $0 \le i <j \le r$. Hence the rational map $\varphi'_1: \rmap \P^{2(r+1)-1}.\P^{\binom{r+1}{2}-1}.$ is defined as $\varphi'_1 ([\mu_{00}:\mu_{01}:\cdots: \mu_{r0}:\mu_{r1}])=[|\begin{smallmatrix}\mu_{00} &\mu_{01}\cr \mu_{10}&\mu_{11}\cr
\end{smallmatrix}|: \cdots: |\begin{smallmatrix}\mu_{r-1,0} &\mu_{r-1,1}\cr \mu_{r0}&\mu_{r1}\cr
\end{smallmatrix}|]$. Denote $\tau_{ij}$ as the coordinates in $\P^{\binom{r+1}{2}-1}$.

By lemma $2.4$, $R_1$ is generated by all $2 \times 2$ determinants $|\begin{smallmatrix}s_{im} &s_{in}\cr s_{jm}&s_{jn}\cr \end{smallmatrix}|$. Blowing up $\P^{(d+1)(r+1)-1}$ along $R_1$, we get $\Gamma_{\varphi_1}$, where $\varphi_1$ is the rational map $\varphi_1: \rmap \P^{(d+1)(r+1)-1}.\P^{\binom{r+1}{2} \cdot \binom{d+1}{2}-1}.$ given by $\varphi_1 ([s_{00}:\cdots:s_{0d}:\cdots: s_{r0}:\cdots:s_{rd}])=[ |\begin{smallmatrix}s_{00} &s_{01}\cr s_{10}&s_{11}\cr \end{smallmatrix}|: \cdots: |\begin{smallmatrix}s_{im} &s_{in}\cr s_{jm}&s_{jn}\cr \end{smallmatrix}|:$

\noindent $\cdots:|\begin{smallmatrix}s_{r-1,d-1} &s_{r-1,d}\cr s_{r,d-1}&s_{r,d}\cr \end{smallmatrix}|]$. Denote $u_{ij, mn}$ as the coordinates in $\P^{\binom{r+1}{2} \cdot \binom{d+1}{2}-1}$, where $0 \le i<j \le r, 0 \le m< n \le d$. 

Define a morphism $F: \map \Gamma_{\varphi'_1} \times \P^{d-1}. \Gamma_{\varphi_1}.$ as $F([\mu_{00}:\mu_{01}:\cdots:\mu_{r0}:\mu_{r1}]; [\tau_{01}:$
$\cdots: \tau_{r-1,r}]; [\nu_0:\cdots:\nu_{d-1}])=[\mu_{0,0} \nu_0: \mu_{0,0} \nu_1+\mu_{0,1} \nu_0:\cdots:\mu_{0,0} \nu_{d-1}+\mu_{0,1} \nu_{d-2}: \mu_{0,1}\nu_{d-1}: \cdots: \mu_{r,0} \nu_0:$ 
$\mu_{r,0} \nu_1+ 
\mu_{r,1} \nu_0:\cdots:
\mu_{r,0} \nu_{d-1}+\mu_{r,1} \nu_{d-2}:\mu_{r,1}\nu_{d-1}]; [\nu_0^2 \tau_{01}: \cdots: \nu_0 \nu_{d-1} \tau_{01}: \cdots: (\nu_m \nu_{n-1}-\nu_{m-1} \nu_n)$
$\tau_{ij}: \cdots: \nu_0 \nu_{d-1} \tau_{r-1,r}:\nu_1 \nu_{d-1} \tau_{r-1,r}:\cdots:\nu_{d-1}^2 \tau_{r-1,r}]$.  It is easy to see that the diagram

$$
\begin{diagram}
\Gamma_{\varphi'_1} \times \P^{d-1} & \rTo^F & \widetilde{R}_2 \subset \Gamma_{\varphi_1}                \\
\dTo^{\pi'_1} &            & \dTo_{\pi_1} \\
\P^{2(r+1)-1} \times \P^{d-1} & \rTo^{\Phi_2} & \P^{(d+1)(r+1)-1}      \\
\end{diagram}
$$
 is commutative.

The image of $F$ included in $\widetilde{R}_2$ can be obtained as follows. Let $E'_1$ be the exceptional divisor of the blow up $\pi'_1: \map \Gamma_{\varphi'_1}.\P^{2(r+1)-1}.$. From
$F(\Gamma_{\varphi'_1} \setminus E'_1 \times \P^{d-1})= \pi_1^{-1} (R_2 \setminus R_1)$, we have that $F(\Gamma_{\varphi'_1} \setminus E'_1 \times \P^{d-1}) \subset \overline{\pi_1^{-1} (R_2 \setminus R_1)}= \widetilde{R}_2$, that is, $\Gamma_{\varphi'_1} \setminus E'_1 \times \P^{d-1} \subset F^{-1} (\widetilde{R}_2)$. Since $F$ is continuous, $F^{-1} (\widetilde{R}_2)$ is closed. Thus $\Gamma_{\varphi'_1} \times \P^{d-1}=\overline{\Gamma_{\varphi'_1} \setminus E'_1 \times \P^{d-1}} \subset F^{-1} (\widetilde{R}_2)$, i.e., $F(\Gamma_{\varphi'_1} \times \P^{d-1}) \subset \widetilde{R}_2$.

Let us show that $F$ is injective.

Suppose two points $P=([\mu_{00}:\mu_{01}:\cdots:\mu_{r0}:\mu_{r1}];[\tau_{01}:\cdots: \tau_{r-1,r}];[\nu_0:\cdots:\nu_{d-1}])$ and $\bar{P}=([\bar{\mu}_{00}:\bar{\mu}_{01}:\cdots:\bar{\mu}_{r0}:\bar{\mu}_{r1}];[\bar{\tau}_{01}:\cdots:\bar{\tau}_{r-1,r}];$
$[\bar{\nu}_0:\cdots:\bar{\nu}_{d-1}])$ have the same image under $F$, we need to show that $P=\bar{P}$. 

From the definition of $F$, $F(P)=F(\bar{P})$ implies that $[\nu_0^2 \tau_{01}: \cdots: \nu_0 \nu_{d-1} \tau_{01}: \cdots: (\nu_m \nu_{n-1}-\nu_{m-1} \nu_n)$
$\tau_{ij}: \cdots: \nu_0 \nu_{d-1} \tau_{r-1,r}:\nu_1 \nu_{d-1} \tau_{r-1,r}:\cdots:\nu_{d-1}^2 \tau_{r-1,r}]=[\bar{\nu}_0^2 \bar{\tau}_{01}: \cdots: \bar{\nu}_0 \bar{\nu}_{d-1} \bar{\tau}_{01}: \cdots: (\bar{\nu}_m \bar{\nu}_{n-1}-\bar{\nu}_{m-1} \bar{\nu}_n)$
$\bar{\tau}_{ij}: \cdots: \bar{\nu}_0 \bar{\nu}_{d-1} \bar{\tau}_{r-1,r}:\bar{\nu}_1 \bar{\nu}_{d-1} \bar{\tau}_{r-1,r}:\cdots:\bar{\nu}_{d-1}^2 \bar{\tau}_{r-1,r}]$.  Let $\tau_{i^{\circ}j^{\circ}}$ be the last non-zero element in the lexicographical order among the coordinates of $P$ in $\P^{\binom{r+1}{2}-1}$ and $\nu_q$ be the last non-zero element in the natural order $\nu_0,\cdots,\nu_{d-1}$ in the coordinates of $P$ in $\P^{d-1}$. Then pick $m=q, n=q+1$, we have $(\nu_m \nu_{n-1}-\nu_{m-1} \nu_n) \tau_{ij}=\nu_q^2 \tau_{i^{\circ}j^{\circ}} \ne 0$. So the corresponding coordinate $(\bar{\nu}_q \bar{\nu}_{(q+1)-1}-\bar{\nu}_{q-1} \bar{\nu}_{q+1}) \bar{\tau}_{i^{\circ} j^{\circ}}=\bar{\nu}_q^2 \bar{\tau}_{i^{\circ}j^{\circ}}$ in $F(\bar{P})$ is also non-zero, which implies that the last non-zero element in the lexicographical order among the coordinates of $\bar{P}$ in $\P^{\binom{r+1}{2}-1}$ is $\bar{\tau}_{i^{\circ}j^{\circ}}$ or an element behind it and the last non-zero element in the natural order $\bar{\nu}_0,\cdots,\bar{\nu}_{d-1}$ in the coordinates of $\bar{P}$ in $\P^{d-1}$ is $\bar{\nu}_q$ or an element behind it. But $P$ and $\bar{P}$ take symmetric roles above, so $\bar{\tau}_{i^{\circ}j^{\circ}}$ is exactly the last non-zero element in the lexicographical order among the coordinates of $\bar{P}$ in $\P^{\binom{r+1}{2}-1}$ and $\bar{\nu}_q$ is the last non-zero element in the natural order $\bar{\nu}_0,\cdots,\bar{\nu}_{d-1}$ in the coordinates of $\bar{P}$ in $\P^{d-1}$. WLOG, take $\nu_q=\bar{\nu}_q=\tau_{i^{\circ}j^{\circ}}=\bar{\tau}_{i^{\circ}j^{\circ}}=1$. From $[\nu_0^2 \tau_{01}: \cdots: \nu_0 \nu_{d-1} \tau_{01}: \cdots: (\nu_m \nu_{n-1}-\nu_{m-1} \nu_n)$
$\tau_{ij}: \cdots: \nu_0 \nu_{d-1} \tau_{r-1,r}:\nu_1 \nu_{d-1} \tau_{r-1,r}:\cdots:\nu_{d-1}^2 \tau_{r-1,r}]=[\bar{\nu}_0^2 \bar{\tau}_{01}: \cdots: \bar{\nu}_0 \bar{\nu}_{d-1} \bar{\tau}_{01}: \cdots: (\bar{\nu}_m \bar{\nu}_{n-1}-\bar{\nu}_{m-1} \bar{\nu}_n)$
$\bar{\tau}_{ij}: \cdots: \bar{\nu}_0 \bar{\nu}_{d-1} \bar{\tau}_{r-1,r}:\bar{\nu}_1 \bar{\nu}_{d-1} \bar{\tau}_{r-1,r}:\cdots:\bar{\nu}_{d-1}^2 \bar{\tau}_{r-1,r}]$, we have that $(\nu_m \nu_{n-1}-\nu_{m-1} \nu_n)\tau_{ij}= (\bar{\nu}_m \bar{\nu}_{n-1}-\bar{\nu}_{m-1} \bar{\nu}_n) \bar{\tau}_{ij}$ for $0 \le i<j \le r, 0 \le m< n \le d$, where $\nu_{-1}=\bar{\nu}_{-1}=\nu_d=\bar{\nu}_d=0$. 
Take $i=i^{\circ}, j=j^{\circ}$ in the above equality, we have $\nu_m \nu_{n-1}-\nu_{m-1} \nu_n= \bar{\nu}_m \bar{\nu}_{n-1}-\bar{\nu}_{m-1} \bar{\nu}_n$ for $0 \le m< n \le d$. And take $m=q, n=q+1$, we have $\tau_{ij}=\bar{\tau}_{ij}$ for $0 \le i<j \le r$.

We claim that $\nu_m \nu_{n-1}-\nu_{m-1} \nu_n$  $(0 \le m< n \le d)$ generate all the monomials of degree $2$ in variables $\nu_0,\cdots,\nu_{d-1}$.

Take $m=0,n=1$, we have $\nu_0 \nu_{1-1}-\nu_{0-1} \nu_1= \nu_0^2$. Similarly, take $m=0, n=j (1<j \le d)$, we have $\nu_0 \nu_{j-1}-\nu_{0-1} \nu_j= \nu_0 \nu_{j-1}$. So $\nu_m \nu_{n-1}-\nu_{m-1} \nu_n$ ($0 \le m< n \le d$) generate all $\nu_0 \nu_i$ for $i=0,\cdots,d-1$.

Suppose $\nu_m \nu_{n-1}-\nu_{m-1} \nu_n$ ($0 \le m< n \le d$) have generated all the monomials in front of $\nu_i \nu_j$ in the lexicographical order for some $0 \le i< j \le d$. Then, take $m=i, n=j+1$, we have that $\nu_i \nu_j= (\nu_i \nu_j-\nu_{i-1} \nu_{j+1})+ \nu_{i-1} \nu_{j+1}$. By our assumption, $\nu_{i-1} \nu_{j+1}$ is generated by $\nu_m \nu_{n-1}-\nu_{m-1} \nu_n$ ($0 \le m< n \le d$), so is $\nu_i \nu_j$. Hence $\nu_m \nu_{n-1}-\nu_{m-1} \nu_n$ ($0 \le m< n \le d$) generate all the monomials of degree $2$ in variables $\nu_0,\cdots,\nu_{d-1}$.

Because all the monomials of degree $2$ in variables $\nu_0,\cdots,\nu_{d-1}$ give the $2$-uple embedding of $\P^{d-1}$ into $\P^{\binom{d+1}{2}-1}$, two points in 
$\P^{d-1}$ will have to coincide if they have the same image in $\P^{\binom{d+1}{2}-1}$ under the $2$-uple embedding.

Now $F(P)=F(\bar{P})$ implies that $\nu_m \nu_{n-1}-\nu_{m-1} \nu_n= \bar{\nu}_m \bar{\nu}_{n-1}-\bar{\nu}_{m-1} \bar{\nu}_n$ for $0 \le m< n \le d$. The previous equation, together with the fact that $\nu_m \nu_{n-1}-\nu_{m-1} \nu_n$ ($0 \le m< n \le d$) generate all the monomials of degree $2$ in variables $\nu_0,\cdots,\nu_{d-1}$, gives that $\nu_i \nu_j=\bar{\nu}_i \bar{\nu}_j$ for all $0 \le i<j \le d-1$. Therefore $\nu_i=\bar{\nu}_i$ for $i=0,\cdots,d-1$. 

By the definition of $\Phi_2$, we see that $[\mu_{00}:\mu_{01}:\cdots:\mu_{r0}:\mu_{r1}]= [\bar{\mu}_{00}:\bar{\mu}_{01}:\cdots:\bar{\mu}_{r0}:\bar{\mu}_{r1}]$ follows from $\pi_1(F(P))=\pi_1(F(\bar{P}))$ and $[\nu_0:\cdots:\nu_{d-1}]=[\bar{\nu}_0:\cdots:\bar{\nu}_d)]$. Hence $F$ is injective.

Let us show that $F$ is onto $\widetilde{R}_2$.

Since $F$ is an injective morphism, $F(\Gamma_{\varphi'_1} \times \P^{d-1})$ is closed in $\Gamma_{\varphi_1}$. From $F(\Gamma_{\varphi'_1} \setminus E'_1 \times \P^{d-1})= \pi_1^{-1} (R_2 \setminus R_1)$, we see that
$\pi_1^{-1} (R_2 \setminus R_1) \subset F(\Gamma_{\varphi'_1} \times \P^{d-1})$. Hence $\widetilde{R}_2$, the closure of $\pi_1^{-1} (R_2 \setminus R_1)$ in $\Gamma_{\varphi_1}$ is a subset of $F(\Gamma_{\varphi'_1} \times \P^{d-1})$. Therefore $F$ is onto $\widetilde{R}_2$.                      

Now $F$ is a bijection between $\Gamma_{\varphi'_1} \times \P^{d-1}$ and $\widetilde{R}_2$. In order to show that $F$ is an isomorphism, it is enough to show that $F$ is a closed immersion.

By a local criterion on closed immersion (see, for example, Proposition 7.3, Chapter II in Hartshorne $[3]$), we need to verify that the coordinate functions in $F$ separate points and tangent vectors. As showed above, the coordinate functions in $F$ separate points. Now let us prove that they also separate tangent vectors.

Since $F$ is a homeomorphism onto $\tilde{R}_2$, to show that the coordinate functions in $F$ separate tangent vectors, it is enough to show that the morphsim of sheaves $\map \CO_{\widetilde{R}_2}. F_{*} \CO_{\Gamma_{\varphi'_1} \times \P^{d-1}}.$ is surjective. We check this surjectivity on stalks.

From the proof of injectivity, we see that $F$ restricted to $\P^{\binom{r+1}{2}-1} \times \P^{d-1}$ gives a morphism $G: \map \P^{\binom{r+1}{2}-1} \times \P^{d-1}. \P^{\binom{r+1}{2} \cdot \binom{d+1}{2}-1}.$, that is,       
$G([\tau_{01}:\cdots:\tau_{r-1,r}];[\nu_0:\cdots:\nu_{d-1}])= [\nu_0^2 \tau_{01}: \cdots: \nu_0 \nu_{d-1} \tau_{01}: \cdots: (\nu_m \nu_{n-1}-\nu_{m-1} \nu_n) \tau_{ij}: \cdots: \nu_0 \nu_{d-1} \tau_{r-1,r}:$
$\nu_1 \nu_{d-1} \tau_{r-1,r}:\cdots:\nu_{d-1}^2 \tau_{r-1,r}]$.  Let $\tilde{G}: \map \P^{\binom{r+1}{2}-1} \times \P^{\binom{d+1}{2}-1}. \P^{\binom{r+1}{2} \cdot \binom{d+1}{2}-1}.$ be the standard Serge embedding. Then $G$ can be decomposed as $\tilde{G} \circ H$, where $H: \map \P^{\binom{r+1}{2}-1} \times \P^{d-1}. \P^{\binom{r+1}{2}-1} \times \P^{\binom{d+1}{2}-1}.$ is the morphism sending $([\tau_{01}:\cdots:\tau_{r-1,r}];[\nu_0:\cdots:\nu_{d-1}])$ into $([\tau_{01}:\cdots:\tau_{r-1,r}];[\nu_0^2:\cdots: (\nu_m \nu_{n-1}-\nu_{m-1} \nu_n):\cdots:\nu_0 \nu_{d-1}]$. Since $\tilde{G}$ is still an embedding when restricted to the image of $H$, so locally $\tau_{ij} ( 0 \le i<j \le r)$ and $\nu_m \nu_{n-1}-\nu_{m-1} \nu_n (0 \le m<n \le d)$ are regular functions in term of the coordinates $u_{ij, mn}$ in $\P^{\binom{r+1}{2} \cdot \binom{d+1}{2}-1}$. Since $\nu_m \nu_{n-1}-\nu_{m-1} \nu_n, 0 \le m<n \le d$ generate all the monomials of degree $2$ in variables $\nu_0,\cdots,\nu_{d-1}$, it follows that any $\nu_i \nu_j$ $(0 \le i \le j \le d-1)$ is locally a regular functions in terms of the coordinates $u_{ij, mn}$ in $\P^{\binom{r+1}{2} \cdot \binom{d+1}{2}-1}$. All the monomials of degree $2$ in variables $\nu_0,\cdots,\nu_{d-1}$ define a $2$-uple embedding of $\P^{d-1}$ into $\P^{\binom{d+1}{2}-1}$, which implies that $\nu_i$ $(i=0,\cdots,d-1)$ are locally regular functions in terms of the coordinates $u_{ij, mn}$ in $\P^{\binom{r+1}{2} \cdot \binom{d+1}{2}-1}$ as well.

Given any point $Q$ in $\widetilde{R}_2$, let $F^{-1}(Q)=P$. We need to show that $\map \CO_{\widetilde{R}_2, Q}. \CO_{\Gamma_{\varphi'_1} \times \P^{d-1}, P}.$ is surjective.

Let $\mu_{i'j'}$, $\tau_{s't'}$ and $\nu_{q'}$ be the last non-zero elements in the lexicographical order among $[\mu_{00}:\mu_{01}:\cdots:\mu_{r0}:\mu_{r1}]$, $[\tau_{01}:\cdots:\tau_{r-1,r}]$ and $[\nu_0:\cdots:\nu_{d-1}]$ in the coordinate of $P$, respectively. Then $\mu_{i'j'}=\tau_{s't'}=\nu_{q'}=1$ gives an affine neighborhood $U$ of $P$. From the definition of $F$, we see that $s_{i',q'+j'}=1$ gives an affine neighborhood $V$ of $Q$. 
Obviously $P$ is also in the open set $1-\nu_{q'-} \nu_{q'+1} \ne 0$. Shrinking the affine neighborhoods if necessary, we may assume that $F(U)=V$ and $1-\nu_{q'-1} \nu_{q'+1} \ne 0$.

Now on $U$ and $V$ we have $\mu_{i0} \nu_{q'} + \nu_{i1} \nu_{q'-1}=s_{iq'}$,  $\mu_{i0} \nu_{q'+1} + \nu_{i1} \nu_{q'}=s_{i,q'+1}$. Since $\nu_{q'}=1$,  it follows from the previous two equations that $\mu_{i0}=\frac{s_{iq'}-\nu_{q'-1}s_{i,q'+1}}{1-\nu_{q'-1} \nu_{q'+1}}$ and $\mu_{i1}=\frac{s_{i,q'+1}-\nu_{q'+1}s_{iq'}}{1-\nu_{q'-1} \nu_{q'+1}}$. Hence $\mu_{i0}, \mu_{i1}$ ($i=0,\cdots,r$) are locally regular functions in terms of $s_{ij}$ ($0 \le i \le r, 0 \le j \le d$) and $\nu_i$ ($i=0, \cdots, d-1$). We have proved that $\nu_i$ ($i=0,\cdots,d-1$) are locally regular functions in terms of the coordinates $u_{ij, mn}$ in $\P^{\binom{r+1}{2} \cdot \binom{d+1}{2}-1}$. So
$\mu_{i0}, \mu_{i1}, i=0,\cdots,r$ are locally regular functions in terms of $s_{ij}$ $(0 \le i \le r, 0 \le j \le d)$ and $u_{ij, mn}$s. Clearly $\mu_{i0}, \mu_{i1}$ $(i=0,\cdots,r)$, $\tau_{ij}$ $(0 \le i <j \le r)$ and $\nu_i$ $(i=0,\cdots,d-1)$ span $\CO_{\Gamma_{\varphi'_1} \times \P^{d-1}, P}$. From all the above, the surjectivity of $\map \CO_{\widetilde{R}_2, Q}. \CO_{\Gamma_{\varphi'_1} \times \P^{d-1}, P}.$ follows. Therefore $F$ is a closed immersion.

That $F$ is a bijection and closed immersion implies $F: \map \Gamma_{\varphi'_1} \times \P^{d-1}. \widetilde{R}_2.$ is an isomorphism (see, for example, J. Harris $[2]$ Corollary $14.10$), thus the smoothness of $\widetilde{R}_2$ follows.

Before we proceed, let us spend some time on setting up notations and explain why we need such ugly ones.

Blowing up $\Gamma_{\varphi_1}$ along $\widetilde{R}_2$, we need a set of generators of $I(\widetilde{R}_2)$. 

Recall that $\varphi_1 ([s_{00}:\cdots:s_{0d}: \cdots: s_{r0}:\cdots: s_{rd}])=[|\begin{smallmatrix}s_{00} &s_{01}\cr s_{10}&s_{11}\cr \end{smallmatrix}|: \cdots: |\begin{smallmatrix}s_{im} &s_{in}\cr s_{jm}&s_{jn}\cr \end{smallmatrix}|: \cdots:$

\noindent $|\begin{smallmatrix}s_{r-1,d-1} &s_{r-1,d}\cr s_{r,d-1}&s_{r,d}\cr \end{smallmatrix}|]$, which sends points outside $R_1$ in $\P^{(d+1)(r+1)-1}$ into $\P^{\binom{r+1}{2} \cdot \binom{d+1}{2}-1}$. We use $u_{ij, mn}$ as the coordinates in $\P^{\binom{r+1}{2} \cdot \binom{d+1}{2}-1}$, where $0 \le i<j \le r, 0 \le m< n \le d$. Over the general points in the image of $\varphi_1$, $[u_{01,01}:\cdots: u_{ij,mn}:\cdots:u_{r-1,r;d-1,d}]=[|\begin{smallmatrix}s_{00} &s_{0,1}\cr s_{10}&s_{11}\cr \end{smallmatrix}|: \cdots: |\begin{smallmatrix}s_{im} &s_{in}\cr s_{jm}&s_{jn}\cr \end{smallmatrix}|: \cdots: |\begin{smallmatrix}s_{r-1,d-1} &s_{r-1,d}\cr s_{r,d-1}&s_{r,d}\cr \end{smallmatrix}|]$.

By Lemma $2.4$, $R_2$ is generated by all the $4 \times 4$ minor determinants of matrices with the form $B_{4,d+2}$. Expanding all such $4 \times 4$ determinants into $2 \times 2$ multiplying $2 \times 2$ determinants, we conclude that the closure $\overline{\varphi_1 (R_2 \setminus R_1)}$ of $\varphi_1 (R_2 \setminus R_1)$ in $\P^{\binom{r+1}{2} \cdot \binom{d+1}{2}-1}$ is generated by $u_{i_1,i_2;m_1,m_2} u_{i_3,i_4; m_3-1, m_4-1}-u_{i_1,i_2;m_1,m_3} u_{i_3,i_4; m_2-1, m_4-1}+u_{i_1,i_2;m_1,m_4} u_{i_3,i_4; m_2-1, m_3-1}+u_{i_1,i_2;m_2,m_3} u_{i_3,i_4; m_1-1, m_4-1}-u_{i_1,i_2;m_2,m_4} u_{i_3,i_4; m_1-1, m_3-1}+u_{i_1,i_2;m_3,m_4} u_{i_3,i_4; m_1-1, m_2-1}$, where $0 \le i_1<i_2 \le r$ and $0 \le i_3<i_4 \le r$ and $0 \le m_1<m_2<m_3<m_4 \le d+1$. The pull-back of $\overline{\varphi_1 (R_2 \setminus R_1)}$ into $\Gamma_{\varphi_1}$ by the projection $\map \Gamma_{\varphi_1}. \P^{\binom{r+1}{2} \cdot \binom{d+1}{2}-1}.$ is $\widetilde{R}_2$. It is easy to see that $\widetilde{R}_2$ is cut out by the same polynomials in $\Gamma_{\varphi_1}$.

The blow up of  $\Gamma_{\varphi_1}$ along $\widetilde{R}_2$ is the graph $\Gamma_{\varphi_2}$ of the rational map $\varphi_2: \Gamma_{\varphi_1} \longrightarrow $

\noindent $\P^{\binom{r+1}{2}^2 \cdot \binom{d+2}{4}-1}$, where $\varphi_2$ maps general points in $\Gamma_{\varphi_1}$ to $[u_{01;01} u_{01;12}-u_{01;02} u_{01; 02}+$

\noindent $u_{01;03} u_{01; 01}: \cdots:u_{i_1,i_2;m_1,m_2} u_{i_3,i_4; m_3-1, m_4-1}-u_{i_1,i_2;m_1,m_3} u_{i_3,i_4; m_2-1, m_4-1}+u_{i_1,i_2;m_1,m_4}$

\noindent $u_{i_3,i_4; m_2-1, m_3-1}+u_{i_1,i_2;m_2,m_3}$
$u_{i_3,i_4; m_1-1, m_4-1}-u_{i_1,i_2;m_2,m_4} u_{i_3,i_4; m_1-1, m_3-1}+u_{i_1,i_2;m_3,m_4}$

\noindent $u_{i_3,i_4; m_1-1, m_2-1}:\cdots:u_{r-1,r;d-2,d-1}$
$u_{r-1,r; d-1, d}-u_{r-1,r;d-2,d} u_{r-1,r; d-2, d} +u_{r-1,r;d-1,d}$

\noindent $u_{r-1,r; d-3, d}]$. We use $u_{i_1,i_2,i_3,i_4; m_1,m_2, m_3, m_4}$ as the coordinates in $\P^{\binom{r+1}{2}^2 \cdot \binom{d+2}{4}-1}$. From the construction of $\varphi_1$ and $\varphi_2$, we see that the coordinate of general point in the image of $\varphi_2$ is equal to $[|M_{01,01;0123}|: \cdots: |M_{i_1i_2,i_3i_4; m_1m_2m_3m_4}|: \cdots: |M_{r-1,r,r-1,r;d-2,d-1,d,d+1}|]$, where $|M_{i_1i_2,i_3i_4; m_1m_2m_3m_4}|$ is the determinant of the $4 \times 4$ submatrix consisting of $(i_1+1)$-th,$(i_2+1)$-th rows from the first block and $(i_3+1)$-th,$(i_4+1)$-th rows from the second block and $(m_j+1)$-th columns ($j=1,\cdots,4$) in $A_{2(r+1) \times (d+2)}$.

By Lemma $2.4$, $R_3$ is generated by $6 \times 6$ minor determinants of submatrices with the form $B_{6,d+3}$. These $6 \times 6$ determinants cannot be expanded as linear combinations of product of $4 \times 4$ determinants of form $|M_{i_1i_2,i_3i_4; m_1m_2m_3m_4}|$, so the defining equations of 

\noindent $\overline{\varphi_2 \circ \pi_1^{-1} (R_3 \setminus R_2)}$ in $\P^{\binom{r+1}{2}^2 \cdot \binom{d+2}{4}-1}$ are not so obvious now. However, we find that we can add a $2 \times 2$ determinant at the corner of the $6 \times 6$ determinant to form a $8 \times 8$ determinant and then expand this $8 \times 8$ determinant into linear combinations of product of two $4 \times 4$ determinants of form $|M_{i_1i_2,i_3i_4; m_1m_2m_3m_4}|$. Specifically we have the following

\begin{equation}
\begin{split} 
\begin{vmatrix} 
s_{i_1,m_1} &\cdots& s_{i_1,m_6}&0&0\cr
s_{i_2,m_1} &\cdots& s_{i_2,m_6}&0&0\cr
s_{i_3,m_1-1} &\cdots& s_{i_3,m_6-1}&0&0\cr
s_{i_4,m_1-1} &\cdots& s_{i_4,m_6-1}&0&0\cr
s_{i_5,m_1-2} &\cdots& s_{i_5,m_6-2}&0&0\cr
s_{i_6,m_1-2} &\cdots& s_{i_6,m_6-2}&0&0\cr
0 &\cdots&0&s_{i_3,n_1}&s_{i_3,n_2}\cr
0 &\cdots& 0&s_{i_4,n_1}&s_{i_4,n_2}\cr
\end{vmatrix} =\begin{vmatrix} 
s_{i_1,m_1} &\cdots& s_{i_1,m_6}&0&0\cr
s_{i_2,m_1} &\cdots& s_{i_2,m_6}&0&0\cr
s_{i_3,m_1-1} &\cdots& s_{i_3,m_6-1}&0&0\cr
s_{i_4,m_1-1} &\cdots& s_{i_4,m_6-1}&0&0\cr
s_{i_5,m_1-2} &\cdots& s_{i_5,m_6-2}&s_{i_5,n_1-1}&s_{i_5,n_2-1}\cr
s_{i_6,m_1-2} &\cdots& s_{i_6,m_6-2}&s_{i_6,n_1-1}&s_{i_6,n_2-1}\cr
s_{i_3,m_1-1} &\cdots& s_{i_3,m_6-1}&s_{i_3,n_1}&s_{i_3,n_2}\cr
s_{i_4,m_1-1} &\cdots& s_{i_4,m_6-1}&s_{i_4,n_1}&s_{i_4,n_2}\cr
\end{vmatrix} 
\end{split}
\end{equation}

Obviously $R_3$ is contained in the common zero locus of all the possible $8 \times 8$ determinant with the above form. The inverse is also true, because if all the $2 \times 2$ determinants at the corner vanish, then the point is in $R_1$, so is in $R_3$; otherwise, at least one of the $2 \times 2$ determinants does not vanish, which implies that all the $6 \times 6$ determinats vanish, so the point is again in $R_3$.

Expanding the $8 \times 8$ determinant on the right hand side of the equation $(3.1)$ into linear combinations of product of two $4 \times 4$ determinants and comparing with the definition of $\varphi_2$, we have that
$\overline{\varphi_2 \circ \pi_1^{-1} (R_3 \setminus R_2)}$ in $\P^{\binom{r+1}{2}^2 \cdot \binom{d+2}{4}-1}$ is generated by $\sum_{1 \le j_1<j_2 \le 6}  (-1)^{j_1+j_2-1}$
$u_{i_1i_2,i_3i_4; m_1,\cdots,\widehat{m_{j_1}},\widehat{m_{j_2}},\cdots,m_6}   u_{i_3i_4,i_5i_6; n_1,n_2,m_{j_1}-1,m_{j_2}-1}$, where $\widehat{x}$ means omitting $x$. It is not hard to see that $\widetilde{R}_3$ in $\Gamma_{\varphi_2}$ is the pull-back of $\overline{\varphi_2 \circ \pi_1^{-1} (R_3 \setminus R_2)}$ by the projection $\map \Gamma_{\varphi_2}. \P^{\binom{r+1}{2}^2 \cdot \binom{d+2}{4}-1}.$ and it is cut out in $\Gamma_{\varphi_2}$ by the same equations.

Now the graph of rational map $\varphi_3: \map \Gamma_{\varphi_2}. \P^{\binom{r+1}{2}^3 \cdot \binom{d+3}{6}\cdot \binom{d+1}{2}-1}.$ gives the blow up of $\Gamma_{\varphi_2}$ along $\widetilde{R}_3$, where $\varphi_3$ sends general points in $\Gamma_{\varphi_2}$ to $[\cdots:\sum_{1 \le j_1<j_2 \le 6}  (-1)^{j_1+j_2-1}$

\noindent $u_{i_1i_2,i_3i_4; m_1,\cdots,\widehat{m_{j_1}},\widehat{m_{j_2}},\cdots,m_6}   u_{i_3i_4,i_5i_6; n_1,n_2,m_{j_1}-1,m_{j_2}-1}:\cdots]$. We use $u_{i_1i_2i_3i_4i_5i_6;m_1m_2m_3m_4m_5m_6}^{n_1,n_2}$ to denote the coordinates in $\P^{\binom{r+1}{2}^3 \cdot \binom{d+3}{6}\cdot \binom{d+1}{2}-1}$, where the up indices represent for the auxiliary determinant $|M_{i_3i_4;n_1n_2}|$. The general point in the image of $\varphi_3$ is $[\cdots:$

\noindent $|M_{i_1i_2i_3i_4i_5i_6;m_1m_2m_3m_4m_5m_6}| \cdot |M_{i_3i_4;n_1n_2}|: \cdots]$, where $|M_{i_1i_2i_3i_4i_5i_6;m_1m_2m_3m_4m_5m_6}|$ is the determinant of the $6 \times 6$ submatrix consisting of $(i_1+1)$-th,$(i_2+1)$-th rows from the first block and $(i_3+1)$-th,$(i_4+1)$-th rows from the second block and $(i_5+1)$-th,$(i_6+1)$-th rows from the third block and $(m_j+1)$-th columns ($j=1,\cdots,6$) in $A_{3(r+1) \times (d+3)}$.

For any natural number $k$, we have the following identity 
$2^{k+1}-2(k+1)= \overset{k-1}{\underset{j=1}{\sum}}2^{j-1} \cdot 2(k-j)$. 
We will use this identity to construct the blow-up 
$\map \Gamma_{\varphi_{k+1}}. \Gamma_{\varphi_k}.$. Specifically for each $2(k+1) \times 2(k+1)$ minor determinant, we can add $2^{j-1}$ auxiliary $2(k-j) \times 2(k-j)$ determinants to it to form a $2^{k+1} \times 2^{k+1}$ determinant, where $j$ runs over $1, \cdots, k-1$. This process of adding auxiliary determinants can be done as follows. Suppose we have known how to add auxiliary determinants to $2k \times 2k$ minor determinants. Let $|M_{i_1,i_2,\cdots,i_{2k+1}i_{2(k+1)};m_1,m_2,\cdots,m_{2k+1},m_{2(k+1)}}|$ be a $2(k+1) \times 2(k+1)$ determinant, and we can first add a $2(k-1) \times 2(k-1)$ determinant $|M_{i_3,i_4,\cdots,i_{2k-1},i_{2k};n_1,n_2,\cdots,n_{2k-3},n_{2(k-1)}}|$ to it at the corner as following

$$
\begin{vmatrix} 
s_{i_1,m_1} & \cdots& s_{i_1,m_{2(k+1)}}&0&\cdots&0\cr
s_{i_2,m_1} &\cdots& s_{i_2,m_{2(k+1)}}&0& \cdots&0\cr
\cdots& \cdots&\cdots&0&\cdots&0\cr
\cdots& \cdots&\cdots&0&\cdots&0\cr
s_{i_{2k-1},m_1-(k-1)} &\cdots& s_{i_{2k-1},m_{2(k+1)}-(k-1)}&0&\cdots&0\cr
s_{i_{2k},m_1-(k-1)} &\cdots& s_{i_{2k},m_{2(k+1)}-(k-1)}&0&\cdots&0\cr
s_{i_{2k+1},m_1-k} &\cdots& s_{i_{2k+1},m_{2(k+1)}-k}&0&\cdots&0\cr
s_{i_{2(k+1)},m_1-k} &\cdots& s_{i_{2(k+1)},m_{2(k+1)}-k}&0&\cdots&0\cr
0 &\cdots&0&s_{i_3,n_1}&\cdots&s_{i_3,n_{2(k-1)}}\cr
0 &\cdots&0&s_{i_4,n_1}&\cdots&s_{i_4,n_{2(k-1)}}\cr
0& \cdots&0&\cdots&\cdots&\cdots\cr
0& \cdots&0&\cdots&\cdots&\cdots\cr
0 &\cdots&0&s_{i_{2k-1},n_1-(k-2)}&\cdots&s_{i_{2k-1},n_{2(k-1)}-(k-2)}\cr
0 &\cdots&0&s_{i_{2k},n_1-(k-2)}&\cdots&s_{i_{2k},n_{2(k-1)}-(k-2)}\cr
\end{vmatrix}
=$$

$$
\begin{vmatrix} 
s_{i_1,m_1} & \cdots& s_{i_1,m_{2(k+1)}}&0&\cdots&0\cr
s_{i_2,m_1} &\cdots& s_{i_2,m_{2(k+1)}}&0& \cdots&0\cr
\cdots& \cdots&\cdots&0&\cdots&0\cr
\cdots& \cdots&\cdots&0&\cdots&0\cr
s_{i_{2k-1},m_1-(k-1)} &\cdots& s_{i_{2k-1},m_{2(k+1)}-(k-1)}&0&\cdots&0\cr
s_{i_{2k},m_1-(k-1)} &\cdots& s_{i_{2k},m_{2(k+1)}-(k-1)}&0&\cdots&0\cr
s_{i_{2k+1},m_1-k} &\cdots& s_{i_{2k+1},m_{2(k+1)}-k}&s_{i_{2k+1},n_1-(k-1)} &\cdots&s_{i_{2k+1},n_{2(k-1)}-(k-1)} \cr
s_{i_{2(k+1)},m_1-k} &\cdots& s_{i_{2(k+1)},m_{2(k+1)}-k}&s_{i_{2(k+1)},n_1-(k-1)} &\cdots&s_{i_{2(k+1)},n_{2(k-1)}-(k-1)} \cr
s_{i_3,m_1-1} &\cdots&s_{i_3,m_{2(k+1)}-1}&s_{i_3,n_1}&\cdots&s_{i_3,n_{2(k-1)}}\cr
s_{i_4,m_1-1}  &\cdots&s_{i_4,m_{2(k+1)}-1}&s_{i_4,n_1}&\cdots&s_{i_4,n_{2(k-1)}}\cr
\cdots& \cdots&\cdots&\cdots&\cdots&\cdots\cr
\cdots& \cdots&\cdots&\cdots&\cdots&\cdots\cr
s_{i_{2k-1},m_1-(k-1)} &\cdots& s_{i_{2k-1},m_{2(k+1)}-(k-1)}&s_{i_{2k-1},n_1-(k-2)}&\cdots&s_{i_{2k-1},n_{2(k-1)}-(k-2)}\cr
s_{i_{2k},m_1-(k-1)} &\cdots& s_{i_{2k},m_{2(k+1)}-(k-1)}&s_{i_{2k},n_1-(k-2)}&\cdots&s_{i_{2k},n_{2(k-1)}-(k-2)}\cr
\end{vmatrix}
$$

Expanding the above determinant from first $2k$ rows, we have that it is a linear combination of product of two $2k \times 2k$ determinants with the desired form. By our assumption, we know how to add auxiliary determinants to such determinants. Thus, the above expansion exactly tells us how to add the lower rank auxiliary determinants to a $2(k+1) \times 2(k+1)$ determinant. So adding auxiliary determinants with the desired form to a $2(k+1) \times 2(k+1)$ determinant can be done. The rational map $\varphi_{k+1}: \Gamma_{\varphi_k} \longrightarrow $
$\P^{\binom{r+1}{2}^{k+1} \cdot \binom{d+(k+1)}{2(k+1)}\cdot \overset{k-1}{\underset{j=1}
{\prod}}\binom{d+(k-j)}{2(k-j)}^{2^{j-1}} 
-1}$ sends general points in $\Gamma_{\varphi_k}$ into 
$[\cdots: \underset {1 \le j_1<j_2 \le 2(k+1)}{\sum} $

\noindent $(-1)^{j_1+j_2-1} u_{i_1,\cdots,i_{2k}; m_1,\cdots,\widehat{m_{j_1}},\widehat{m_{j_2}},\cdots,m_{2(k+1)}}^{J_k}   u_{i_3,\cdots,i_{2(k+1)}; n_1,\cdots,n_{2(k-1)},m_{j_1}-1,m_{j_2}-1}^{\tilde{J}_k}:$
$\cdots]$, where $J_k$ and $\tilde{J}_k$ are indices representing the auxiliary determinants added to $2k \times 2k$ determinants. The coordinate of the general point in the image of $\varphi_{k+1}$ is given by the product of $|M_{i_1i_2,\cdots,i_{2k+1}i_{2(k+1)};m_1,m_2,\cdots,m_{2k+1}m_{2(k+1)}}|$ and those auxiliary determinants added to it. We use 
$u_{i_1,\cdots,i_{2(k+1)};m_1,\cdots,m_{2(k+1)}}^{n_1,\cdots,n_{2(k-1)},J_k,\tilde{J}_k}=u_{i_1,\cdots,i_{2(k+1)};m_1,\cdots,m_{2(k+1)}}^{J_{k+1}}$ as the coordinates in

\noindent $\P^{\binom{r+1}{2}^{k+1} \cdot \binom{d+(k+1)}{2(k+1)}\cdot \overset{k-1}{\underset{j=1}{\prod}}\binom{d+(k-j)}{2(k-j)}^{2^{j-1}} 
-1}$.

Now let us continue our proof that the iterated blow-up can be carried out.

We have showed that $\widetilde{R}_2$ is smooth, so the second step in the iterated blow-ups can be done. Suppose that the iterated blow-ups can be carried out in the first $k$ $(k \ge 2)$ steps. If $k=d$, we are done. If $k<d$, let us show that we can proceed to the $(k+1)$-th step. Since $\Gamma_{\varphi_{k}}$ is smooth, in order to carry out the $(k+1)$-th blow up, we only need to show that the proper transformation $\widetilde{R}_{k+1}$ of $R_{k+1}$ in $\Gamma_{\varphi_{k}}$ is smooth. We use the same idea as we did before in showing the smoothness of $\widetilde{R}_{2}$. We construct a birational morphism $\Phi_{k+1}: \map \P^{(k+1)(r+1)-1} \times \P^{d-k}. R_{k+1} \subset \P^{(d+1)(r+1)-1}.$ with $\Phi([\mu_{0,0}:\cdots:\mu_{0,k}:\mu_{1,0}:\cdots:\mu_{1,k}:$
$\cdots: \mu_{r,0}:\cdots:\mu_{r,k}]; [\nu_0: \cdots: \nu_{d-k}])=[ \mu_{0,0} \nu_0: \mu_{0,0} \nu_1+\mu_{0,1} \nu_0:\cdots:\underset{p+q=s}{\sum} \mu_{0,p} \nu_{q}: \cdots:\mu_{0,k-1} \nu_{d-k} $
$+\mu_{0,k}\nu_{d-k-1}: \mu_{0,k}\nu_{d-k}: \cdots: \mu_{r,0} \nu_0: \mu_{r,0} \nu_1+\mu_{r,1} \nu_0:\cdots:\underset{p+q=s}{\sum} \mu_{r,p} \nu_{q}: \cdots:\mu_{r,k-1} \nu_{d-k}+\mu_{r,k} \nu_{d-k-1}:$
$\mu_{0,k}\nu_{d-k}]$. The space $\P^{(k+1)(r+1)-1}$ parametrizes $(r+1)$-tuples $(g_0,\cdots,g_r)$ modulo homothety and the morphism $\Phi$ sends $(g_0,\cdots,g_r;f)$ to $(g_0 \cdot f,\cdots,g_r \cdot f)$, where $g_i (x,y)=\underset{p=0}{\overset{k}{\sum}} \mu_{i,p} x^{k-p} y^p (i$
$=0, \cdots, r)$ are homogeneous polynomials of degree $k$. In $\P^{(k+1)(r+1)-1}$, the set parametrizing $(g_0,\cdots,g_r)$ with at least one common root has a natural stratification $R'_1 \subset \cdots \subset R'_{k}$. By our inductive assumption, the iterated blow-ups can be carried out on $\P^{(k+1)(r+1)-1}$. Blow-up $\P^{(k+1)(r+1)-1}$ along $R'_1$ gives $\pi'_1: \map \Gamma_{\varphi'_{1}}. \P^{(k+1)(r+1)-1}.$. Continue the iterated blow-ups until we get $\pi'_{k}: \map \Gamma_{\varphi'_{k}}.\Gamma_{\varphi'_{k-1}}.$. Since at each step we blow up a smooth variety along a smooth center, the final outcome $\Gamma_{\varphi'_{k}}$ is smooth. We will show that we can lift the morphism $\Phi_{k+1}: \map \P^{(k+1)(r+1)-1} \times \P^{d-k}. R_{k+1} \subset \P^{(d+1)(r+1)-1}.$ to a morphism $F: \map \Gamma_{\varphi'_{k}} \times \P^{d-k}. \widetilde{R}_{k+1} \subset \Gamma_{\varphi_{k}}.$ such that the following diagram 

$$
\begin{diagram}
\Gamma_{\varphi'_{k}} \times \P^{d-k} & \rTo^F & \widetilde{R}_{k+1} \subset \Gamma_{\varphi_{k}}\\
\dTo^{\pi'} &            & \dTo_{\pi} \\
\P^{(k+1)(r+1)-1} \times \P^{d-k} & \rTo^{\Phi_{k+1}} & \P^{(d+1)(r+1)-1}      \\
\end{diagram}
$$
 is commutative, where $\pi'=\pi'_1 \circ \pi'_2 \circ \cdots \circ \pi'_{k}$ and $\pi=\pi_1 \circ \pi_2 \circ \cdots \circ \pi_{k}$, respectively. Moreover we will show that
$F$ gives an isomorphim between $\Gamma_{\varphi'_{k}} \times \P^{d-k}$ and $\widetilde{R}_{k+1}$, therefore $\widetilde{R}_{k+1}$ is smooth.

For any point $P$ in $\Gamma_{\varphi'_{k}} \times \P^{d-k} \subset \P^{(k+1)(r+1)-1} \times \P^{\binom{r+1}{2} \cdot \binom{k+1}{2}-1}\times \P^{\binom{r+1}{2}^2 \cdot \binom{k+2}{4}-1} \times$

\noindent $\overset{k}{\underset{s=3}{\prod}}  \P^{\binom{r+1}{2}^{s} \cdot \binom{k+s}{2s}\cdot \overset{k-2}{\underset{j=1}
{\prod}}\binom{k+(s-1-j)}{2(s-1-j)}^{2^{j-1}} 
-1} \times \P^{d-k}$, let us define its image $F(P)$ in $\Gamma_{\varphi_{k}}$. 

Since $\Gamma_{\varphi_{k}} \subset \P^{(d+1)(r+1)-1} \times \P^{\binom{r+1}{2} \cdot \binom{d+1}{2}-1}\times \P^{\binom{r+1}{2}^2 \cdot \binom{d+2}{4}-1} \times \overset{k}{\underset{s=3}{\prod}} $

\noindent $\P^{\binom{r+1}{2}^{s} \cdot \binom{d+s}{2s}\cdot \overset{s-2}{\underset{j=1}
{\prod}}\binom{d+(s-1-j)}{2(s-1-j)}^{2^{j-1}} 
-1}$, we will give the component of $F(P)$ in each projective space.

For convenience, let $\mu_{ij}=0$ for $i=0,\cdots,r, j=k+1,\cdots,d$ and $\nu_i=0$ for $i=d-k+1,\cdots,d$. Then the component of $F(P)$ in $\P^{(d+1)(r+1)-1}$ is $[ \underset{p+q=0}{\sum} \mu_{0,p} \nu_{q}:\underset{p+q=1}{\sum} \mu_{0,p} \nu_{q}:\cdots: \underset{p+q=d}{\sum} \mu_{0,p} \nu_{q}:\cdots: \underset{p+q=0}{\sum} \mu_{r,p} \nu_{q}:   
\underset{p+q=1}{\sum} \mu_{r,p} \nu_{q}:\cdots: \underset{p+q=d}{\sum} \mu_{r,p} \nu_{q}]$, which  is exactly $\Phi_{k+1}$ when we restrict $F$ to $\P^{(k+1)(r+1)-1} \times \P^{d-k}$.
 
The component of $F(P)$ in $\P^{\binom{r+1}{2} \cdot \binom{d+1}{2}-1}$ is given as follows. 

The component in $\P^{\binom{r+1}{2} \cdot \binom{d+1}{2}-1}$ of the coordinate of a general point in $\Gamma_{\varphi_{1}}$ is $[\cdots:u_{i_1i_2,m_1m_2}:\cdots]=[\cdots:|M_{i_1i_2,m_1m_2}|:\cdots]=[\cdots:|\begin{smallmatrix}s_{i_1,m_1} &s_{i_1,m_2}\cr s_{i_2,m_1}&s_{i_2,m_2}\cr \end{smallmatrix}|:\cdots]= [\cdots: $

\noindent $|\begin{smallmatrix} \sum_{p_1+q_1=m_1}  \mu_{i_1,p_1} \nu_{q_1}& \sum_{p_2+q_2=m_2}  \mu_{i_1,p_2} \nu_{q_2}\cr \sum_{p_1+q_1=m_1}  \mu_{i_2,p_1} \nu_{q_1}& \sum_{p_2+q_2=m_2}  \mu_{i_2,p_2} \nu_{q_2}\cr \end{smallmatrix}|:\cdots]=[\cdots: \underset{p_1+q_1=m_1}{\sum} \underset{p_2+q_2=m_2}{\sum} \nu_{q_1} \nu_{q_2} |\begin{smallmatrix} \mu_{i_1,p_1}  & \mu_{i_1,p_2} \cr \mu_{i_2,p_1} & \mu_{i_2,p_2} \cr \end{smallmatrix}|:\cdots]=[\cdots: \underset{p_1+q_1=m_1}{\sum} \underset{p_2+q_2=m_2}{\sum} \nu_{q_1} \nu_{q_2} \tau_{i_1i_2,p_1p_2}:\cdots]  $, so we take the coordinate of $F(P)$ at the position $u_{i_1i_2,m_1m_2}$ as $\underset{p_1+q_1=m_1}{\sum} \underset{p_2+q_2=m_2}{\sum} \nu_{q_1} \nu_{q_2} \tau_{i_1i_2,p_1p_2}$, where $\tau_{i_1i_2,p_1p_2}$ are coordinates in $\P^{\binom{r+1}{2} \cdot \binom{k+1}{2}-1}$. We need to verify that at least one of $\underset{p_1+q_1=m_1}{\sum} \underset{p_2+q_2=m_2}{\sum} \nu_{q_1} \nu_{q_2}$
$\tau_{i_1i_2,p_1p_2} \ne 0$, this can be done easily. For example, there exist $i_1$ and $i_2$ such that $\tau_{i_1i_2,p_1p_2} \ne 0$ for some $p_1< p_2$. Fix $i_1$ and $i_2$ and take lexicographical order among all letters $p_1p_2 (p_1<p_2)$, let $\tau_{i_1i_2,p_1^{\circ} p_2^{\circ}}$ be the last nonzero element. Let $\nu_q$ be the last nonzero element among $\nu_0,\cdots,\nu_{d-k}$. Then let $q+p_1^{\circ}=m_1, q+p_2^{\circ}=m_2$, we have   $\underset{p_1+q_1=m_1}{\sum} \underset{p_2+q_2=m_2}{\sum} \nu_{q_1} \nu_{q_2} \tau_{i_1i_2,p_1p_2}=\nu_q^2 \tau_{i_1i_2,p_1^{\circ} p_2^{\circ}} \ne 0$. 

The component in $\P^{\binom{r+1}{2}^2 \cdot \binom{d+2}{4}-1}$ of the coordinate of a general point in $\Gamma_{\varphi_{2}}$ is $[\cdots:u_{i_1i_2i_3i_4;m_1m_2m_3m_4}:\cdots]=[\cdots:|M_{i_1i_2i_3i_4;m_1m_2m_3m_4}|:\cdots]$, where

$$
|M_{i_1i_2i_3i_4;m_1m_2m_3m_4}|=
\begin{vmatrix} 
s_{i_1,m_1} & s_{i_1,m_2} & s_{i_1,m_3} &s_{i_1,m_4} \cr
s_{i_2,m_1} & s_{i_2,m_2} & s_{i_2,m_3} &s_{i_2,m_4} \cr
s_{i_3,m_1-1} & s_{i_3,m_2-1} & s_{i_3,m_3-1} &s_{i_3,m_4-1} \cr
s_{i_4,m_1-1} & s_{i_4,m_2-1} & s_{i_4,m_3-1} &s_{i_4,m_4-1} \cr
\end{vmatrix}
$$

Since $[\cdots:s_{i,m}:\cdots]=[\cdots: \overset{m}{\underset{q=0}{\sum}} \mu_{i,m-q} \nu_q :\cdots]$, then $[\cdots:|M_{i_1i_2i_3i_4;m_1m_2m_3m_4}|:\cdots]=[$

$$
\cdots:
\begin{vmatrix} 
\overset{m_1}{\underset{q_1=0}{\sum}} \mu_{i_1,m_1-q_1} \nu_{q_1} & \overset{m_2}{\underset{q_2=0}{\sum}} \mu_{i_1,m_2-q_2} \nu_{q_2} & \overset{m_3}{\underset{q_3=0}{\sum}} \mu_{i_1,m_3-q_3} \nu_{q_3} &\overset{m_4}{\underset{q_4=0}{\sum}} \mu_{i_1,m_4-q_4} \nu_{q_4} \cr
\overset{m_1}{\underset{q_1=0}{\sum}} \mu_{i_2,m_1-q_1} \nu_{q_1} & \overset{m_2}{\underset{q_2=0}{\sum}} \mu_{i_2,m_2-q_2} \nu_{q_2} & \overset{m_3}{\underset{q_3=0}{\sum}} \mu_{i_2,m_3-q_3} \nu_{q_3} &\overset{m_4}{\underset{q_4=0}{\sum}} \mu_{i_2,m_4-q_4} \nu_{q_4} \cr
\overset{m_1-1}{\underset{q_1=0}{\sum}} \mu_{i_3,m_1-1-q_1} \nu_{q_1} & \overset{m_2-1}{\underset{q_2=0}{\sum}} \mu_{i_3,m_2-1-q_2} \nu_{q_2} & \overset{m_3-1}{\underset{q_3=0}{\sum}} \mu_{i_3,m_3-1-q_3} \nu_{q_3} &\overset{m_4-1}{\underset{q_4=0}{\sum}} \mu_{i_3,m_4-1-q_4} \nu_{q_4} \cr
\overset{m_1-1}{\underset{q_1=0}{\sum}} \mu_{i_4,m_1-1-q_1} \nu_{q_1} & \overset{m_2-1}{\underset{q_2=0}{\sum}} \mu_{i_4,m_2-1-q_2} \nu_{q_2} & \overset{m_3-1}{\underset{q_3=0}{\sum}} \mu_{i_4,m_3-1-q_3} \nu_{q_3} &\overset{m_4-1}{\underset{q_4=0}{\sum}} \mu_{i_4,m_4-1-q_4} \nu_{q_4} \cr
\end{vmatrix} :\cdots]$$

$$
=[\cdots: \overset{m_1}{\underset{q_1=0}{\sum}} \overset{m_2}{\underset{q_2=0}{\sum}}  \overset{m_3}{\underset{q_3=0}{\sum}} \overset{m_4}{\underset{q_4=0}{\sum}} \nu_{q_1} 
\nu_{q_2}\nu_{q_3} \nu_{q_4} 
\begin{vmatrix} 
\mu_{i_1,m_1-q_1} & \mu_{i_1,m_2-q_2} & \mu_{i_1,m_3-q_3} & \mu_{i_1,m_4-q_4} \cr
\mu_{i_2,m_1-q_1} & \mu_{i_2,m_2-q_2} & \mu_{i_2,m_3-q_3} & \mu_{i_2,m_4-q_4} \cr
\mu_{i_3,m_1-1-q_1} & \mu_{i_3,m_2-1-q_2} & \mu_{i_3,m_3-1-q_3} & \mu_{i_3,m_4-1-q_4} \cr
\mu_{i_4,m_1-1-q_1} & \mu_{i_4,m_2-1-q_2} & \mu_{i_4,m_3-1-q_3} & \mu_{i_4,m_4-1-q_4} \cr
\end{vmatrix}:\cdots
 $$
$]=[\cdots: \overset{m_1}{\underset{q_1=0}{\sum}} \overset{m_2}{\underset{q_2=0}{\sum}}  \overset{m_3}{\underset{q_3=0}{\sum}} \overset{m_4}{\underset{q_4=0}{\sum}} \nu_{q_1} 
\nu_{q_2}\nu_{q_3} \nu_{q_4} \tau_{i_1i_2i_3i_4;m_1-q_1,m_2-q_2,m_3-q_3,m_4-q_4}:\cdots]$. In the above computation we use the property that any determinant vanishes if it contains two repeated columns, and we also use $\tau_{i_1i_2i_3i_4;m_1-q_1,m_2-q_2,m_3-q_3,m_4-q_4}=0$ if any two of $m_i-q_i, i=1,2,3,4$ are equal. So we will take the coordinate of $F(P)$ at the position $u_{i_1i_2i_3i_4;m_1m_2m_3m_4}$ as $\overset{m_1}{\underset{q_1=0}{\sum}} \overset{m_2}{\underset{q_2=0}{\sum}}  \overset{m_3}{\underset{q_3=0}{\sum}} \overset{m_4}{\underset{q_4=0}{\sum}} \nu_{q_1} 
\nu_{q_2}\nu_{q_3} \nu_{q_4}$
$\tau_{i_1i_2i_3i_4;m_1-q_1,m_2-q_2,m_3-q_3,m_4-q_4}$. Use the same method as in  showing at least one of $\underset{p_1+q_1=m_1}{\sum} \underset{p_2+q_2=m_2}{\sum}$
$\nu_{q_1} \nu_{q_2} \tau_{i_1i_2,p_1p_2} \ne 0$, we can prove that at least one of $\overset{m_1}{\underset{q_1=0}{\sum}} \overset{m_2}{\underset{q_2=0}{\sum}}  \overset{m_3}{\underset{q_3=0}{\sum}} \overset{m_4}{\underset{q_4=0}{\sum}} \nu_{q_1} 
\nu_{q_2}\nu_{q_3} \nu_{q_4}$
$\tau_{i_1i_2i_3i_4;m_1-q_1,m_2-q_2,m_3-q_3,m_4-q_4} \ne 0$ as well.

Now let us specify the component of $F(P)$ in each $\P^{\binom{r+1}{2}^{s} \cdot \binom{k+s}{2s}\cdot \overset{k-2}{\underset{j=1}
{\prod}}\binom{k+(s-1-j)}{2(s-1-j)}^{2^{j-1}} 
-1}$ for $s=3, \cdots,k$.

Over the general point of the image $\varphi_{s}: \rmap \Gamma_{\varphi_{s-1}}.\P^{\binom{r+1}{2}^{s} \cdot \binom{k+s}{2s}\cdot \overset{k-2}{\underset{j=1}
{\prod}}\binom{k+(s-1-j)}{2(s-1-j)}^{2^{j-1}} 
-1}.$, the coordinate
$[\cdots: u_{i_1,\cdots,i_{2s};m_1,\cdots,m_{2s}}^{J_{s}}:\cdots]=[\cdots:\underset{|M|}{\prod} |M|:\cdots]$, where $|M|$ runs over $|M_{i_1,\cdots,i_{2s};m_1,\cdots,m_{2s}}|$ and those   
auxiliary determinants added to $|M_{i_1,\cdots,i_{2s};m_1,\cdots,m_{2s}}|$ with indices in $J_s$.
Experienced reader may have seen that $[\cdots:\underset{|M|}{\prod} |M|:\cdots]=[\cdots: \underset{q_i}{\sum}^{(2s)} \underset{\tilde{q}_j}{\sum}^{(2^s-2s)} \prod \nu_{q_i} \cdot \prod \nu_{\tilde{q_j}} \tau_{i_1,\cdots,i_{2s};m_1-q_1,\cdots,m_{2s}-q_{2s}}^{\cdots,n_j-\tilde{q_j},\cdots}:\cdots]$, where $\underset{q_i}{\sum}^{(2s)}= \overset{m_1}{\underset{q_1=0}{\sum}} \cdots \overset{m_{2s}}{\underset{q_{2s}=0}{\sum}}$ and $\underset{\tilde{q}_j}{\sum}^{(2^s-2s)}$ represents $2^s-2s$ sums over all $0 \le \tilde{q}_j \le n_j$, $n_j$ appearing in $J_s$ as the column index $n_j+1$ in $A_{(r+1)p \times (d+p)}$ for the added auxiliary $2p \times 2p$ $(1 \le p \le s-2)$ minor determinants. So $\prod \nu_{q_i} \cdot \prod \nu_{\tilde{q_j}}$ is a monomial in $\nu_0,\cdots,\nu_{d-k}$ of degree $2^s$.

We can order the indices in $u_{i_1,\cdots,i_{2s};m_1,\cdots,m_{2s}}^{J_{s}}$ as follows.

Indices $i_1,\cdots,i_{2s};m_1,\cdots,m_{2s}$ appear in the first position. The group of the column indices of the added $2(s-2) \times 2(s-2)$ determinant appear in the second position. The two groups of the column indices of the two added $2(s-3)\times 2(s-3)$ determinants should be put in the third position but we have to decide which group is put first, and this can be done by comparing their row indices. In general, the groups of the column indices of the $2^{j-1}$ added $2(s-1-j) \times 2(s-1-j)$ determinants appear in the $(j+1)$-th position. The priority among the group of indices in the $(j+1)$-th position is given by their row indices appeared in each $2(s-1-j) \times 2(s-1-j)$ determinant. Now we can arrange $u_{i_1,\cdots,i_{2s};m_1,\cdots,m_{2s}}^{J_{s}}$ in $\P^{\binom{r+1}{2}^{s} \cdot \binom{d+s}{2s}\cdot \overset{d-2}{\underset{j=1}
{\prod}}\binom{d+(s-1-j)}{2(s-1-j)}^{2^{j-1}} 
-1}$ in the lexicographical order. Essentially as we did before, we can show that at least one of  
$\underset{q_i}{\sum}^{(2s)} \underset{\tilde{q}_j}{\sum}^{(2^s-2s)} \prod \nu_{q_i} \cdot \prod \nu_{\tilde{q_j}} \tau_{i_1,\cdots,i_{2s};m_1-q_1,\cdots,m_{2s}-q_{2s}}^{\cdots,n_j-\tilde{q_j},\cdots} \ne 0$ for some $i_1, \cdots, i_{2s}$, $m_1,\cdots,m_{2s}$ and $J_s$. We leave this as an exercise to the interest readers.

From above we see that $F: \map \Gamma_{\varphi'_{k}} \times \P^{d-k}.\Gamma_{\varphi_{k}}.$ is a morphism. We will show that it is in fact an isomorphism from $\Gamma_{\varphi'_{k}} \times \P^{d-k}$ onto $\widetilde{R}_{k+1}$.

First let us show that $F$ is injective.

Suppose that $F(P)=F(\bar{P})$, we show that $P=\bar{P}$.

From the definition of $F$, we see that $[\cdots: \underset{q_i}{\sum}^{(2k)} \underset{\tilde{q}_j}{\sum}^{(2^k-2k)} \prod \nu_{q_i} \cdot \prod \nu_{\tilde{q_j}} $

\noindent $\tau_{i_1,\cdots,i_{2k};m_1-q_1,\cdots,m_{2k}-q_{2k}}^{\cdots,n_j-\tilde{q_j},\cdots}:\cdots]=[\cdots: \underset{q_i}{\sum}^{(2k)} \underset{\tilde{q}_j}{\sum}^{(2^k-2k)} \prod \bar{\nu}_{q_i} \cdot $
$\prod \bar{\nu}_{\tilde{q_j}} \bar{\tau}_{i_1,\cdots,i_{2k};m_1-q_1,\cdots,m_{2k}-q_{2k}}^{\cdots,n_j-\tilde{q_j},\cdots}$

\noindent $:\cdots]$, where $\bar{\nu}$ and $\bar{\tau}$ represent the coordinates appearing in $\bar{P}$. By the construction of $\Gamma_{\varphi'_{k}}$,  the indices $\{m_1-q_1, \cdots, m_{2k}-q_{2k}\}$ of $\tau_{i_1,\cdots,i_{2s};m_1-q_1,\cdots,m_{2k}-q_{2k}}^{\cdots,n_j-\tilde{q_j},\cdots}$ and $\bar{\tau}_{i_1,\cdots,i_{2k};m_1-q_1,\cdots,m_{2k}-q_{2k}}^{\cdots,n_j-\tilde{q_j},\cdots}$ in $\P^{\binom{r+1}{2}^{k} \cdot \binom{k+k}{2k}\cdot \overset{k-2}{\underset{j=1}
{\prod}}\binom{k+(k-1-j)}{2(k-1-j)}^{2^{j-1}} 
-1}= \P^{\binom{r+1}{2}^{k} \cdot \overset{k-2}{\underset{j=1}
{\prod}}\binom{2k-1-j}{2(k-1-j)}^{2^{j-1}} 
-1}$ have only one choice, that is, $\{m_1-q_1,\cdots,m_{2k}-q_{2k}\}=\{0,\cdots,2k-1\}$. We will show that $[\cdots:\tau_{i_1,\cdots,i_{2k};0,\cdots,2k-1}^{J_{k}}:\cdots]=[\cdots:\bar{\tau}_{i_1,\cdots,i_{2k};0,\cdots,2k-1}^{J_{k}}:\cdots]$ and $[\nu_0:\cdots:\nu_{d-k}]=[\bar{\nu}_0:\cdots:\bar{\nu}_{d-k}]$. 

For each fixed $\{m_1,\cdots,m_{2k}\}$ with $0 \le m_1<m_2<\cdots<m_{2k} \le d+k-1$, we have 
$\underset{q_i}{\sum}^{(2k)} \underset{\tilde{q}_j}{\sum}^{(2^k-2k)} \prod \nu_{q_i} \cdot \prod \nu_{\tilde{q_j}} \tau_{i_1,\cdots,i_{2k};m_1-q_1,\cdots,m_{2k}-q_{2k}}^{\cdots,n_j-\tilde{q_j},\cdots}= \underset{\{l_1,\cdots,l_{2k}\}=\{0,\cdots,2k-1\}}{\sum} \underset{\tilde{q}_j}{\sum}^{(2^k-2k)} \prod \nu_{m_i-l_i} \cdot \prod \nu_{\tilde{q_j}} \tau_{i_1,\cdots,i_{2k};l_1,\cdots,l_{2k}}^{\cdots,n_j-\tilde{q_j},\cdots}=$ 
$(\underset{\{l_1,\cdots,l_{2k}\}=
\{0,\cdots,2k-1\}}{\sum} \pm\prod \nu_{m_i-l_i})\cdot (\underset{\tilde{q}_j}{\sum}^{(2^k-2k)} \prod \nu_{\tilde{q_j}} \tau_{i_1,\cdots,i_{2k};0,\cdots,2k-1}^{\cdots,n_j-\tilde{q_j},\cdots})$, where the sign $\pm$ depends on the permutation $l_1 \cdots l_{2k}$ even or odd.

Now let us show that the linear combinations of $\underset{\{l_1,\cdots,l_{2k}\}=
\{0,\cdots,2k-1\}}{\sum} \pm\prod \nu_{m_i-l_i}$ for $0 \le m_1<m_2<\cdots<m_{2k} \le d+k-1$ produce all the monomials of degree $2k$ in $\nu_0,\cdots,\nu_{d-k}$.

Take $m_i=i-1, i=1,\cdots,2k$, then $\underset{\{l_1,\cdots,l_{2k}\}=
\{0,\cdots,2k-1\}}{\sum} \pm\prod \nu_{m_i-l_i}=\nu_0^{2k}$.  Take $m_i=i-1, i=1,\cdots,2k-1$ and $m_{2k}=2k-1+j$ for $1 \le j \le d-k$, then $\underset{\{l_1,\cdots,l_{2k}\}=
\{0,\cdots,2k-1\}}{\sum} \pm\prod \nu_{m_i-l_i}=\nu_0^{2k-1} \nu_j$.  Take $m_i=i-1, i=1,\cdots,2k-2$ and $m_{2k-1}=2k-1, m_{2k}=2k$, then $\underset{\{l_1,\cdots,l_{2k}\}=
\{0,\cdots,2k-1\}}{\sum} \pm\prod \nu_{m_i-l_i}=\nu_0^{2k-2} \nu_1^2-\nu_0^{2k-1} \nu_2$.  Since $\nu_0^{2k-1} \nu_2$ is a linear combination of $\underset{\{l_1,\cdots,l_{2k}\}=
\{0,\cdots,2k-1\}}{\sum} \pm\prod \nu_{m_i-l_i}$, so is $\nu_0^{2k-2} \nu_1^2$. Take $m_i=i-1, i=1,\cdots,2k-2$ and $m_{2k-1}=2k-1, m_{2k}=2k-1+j$ for $2 \le j \le d-k$, then $\underset{\{l_1,\cdots,l_{2k}\}=
\{0,\cdots,2k-1\}}{\sum} \pm\prod \nu_{m_i-l_i}=\nu_0^{2k-2} \nu_1 \nu_{j}-\nu_0^{2k-1} \nu_{j+1}$.  Since $\nu_0^{2k-1} \nu_{j+1}$ is a linear combination of $\underset{\{l_1,\cdots,l_{2k}\}=
\{0,\cdots,2k-1\}}{\sum} \pm\prod \nu_{m_i-l_i}$, so is $\nu_0^{2k-2} \nu_1 \nu_j$.  
 
In general we can arrange all the monomials $\nu_0^{i_0} \cdot \cdots \nu_{d-k}^{i_{d-k}}$ in the lexicographical order, where $i_0+\cdots+i_{d-k}=2k$. Suppose that all the monomials in front of $\nu_0^{i_0} \cdot \cdots \nu_{d-k}^{i_{d-k}}$ are linear combination of $\underset{\{l_1,\cdots,l_{2k}\}=
\{0,\cdots,2k-1\}}{\sum} \pm\prod \nu_{m_i-l_i}$. Let the ordered tuple $(m_1,\cdots,m_{2k})=(0,\cdots,i_0-1;i_0+1,\cdots,i_0+i_1;\cdots,i_0+i_1+\cdots+i_t+t+1,\cdots,i_0+i_1+\cdots+i_{t}+i_{t+1}+t;\cdots, i_0+i_1+\cdots+i_{d-k-1}+d-k,\cdots,i_0+i_1+\cdots+i_{d-k-1}+i_{d-k}+d-k-1)$. Note that some piece $i_0+i_1+\cdots+i_t+t+1,\cdots,i_0+i_1+\cdots+i_{t}+i_{t+1}+t$ may not actually appear in $(m_1,\cdots,m_{2k})$ if $i_{t+1}=0$, so $m_{2k}$ may be smaller than $i_0+i_1+\cdots+i_{d-k-1}+i_{d-k}+d-k-1=2k+d-k-1=d+k-1$. It is easy to see that $0 \le m_1<m_2<\cdots<m_{2k} \le d+k-1$. For such a $(m_1,\cdots,m_{2k})$, we have that $\nu_0^{i_0} \cdot \cdots \nu_{d-k}^{i_{d-k}}$ +  $\underset{\{l_1,\cdots,l_{2k}\}=
\{0,\cdots,2k-1\}}{\sum'} \pm\prod \nu_{m_i-l_i}=0$, where $\sum'$ means that the sum is over all $\{l_1,\cdots,l_{2k}\}=\{0,\cdots,2k-1\}$ with at least one of $l_i \ne i-1$. Pick any monomial in $\sum'$ and let $l_t$ be the first one in $(l_1,\cdots,l_{2k})$ with the property $l_t \ne t-1$, then $l_t > t-1$. Thus $m_t-l_t<m_t-(t-1)$. So a smaller factor $\nu_{m_t-l_t}$ appears in the monomial, which implies that this monomial is in front of $\nu_0^{i_0} \cdot \cdots \nu_{d-k}^{i_{d-k}}$. By our inductive assumption, the monomials in $\sum'$ are linear combinations of $\underset{\{l_1,\cdots,l_{2k}\}=
\{0,\cdots,2k-1\}}{\sum} \pm\prod \nu_{m_i-l_i}$, so is $\nu_0^{i_0} \cdot \cdots \nu_{d-k}^{i_{d-k}}$.

Now from $[\cdots:(\underset{\{l_1,\cdots,l_{2k}\}=
\{0,\cdots,2k-1\}}{\sum} \pm\prod \nu_{m_i-l_i})\cdot (\underset{\tilde{q}_j}{\sum}^{(2^k-2k)} \prod \nu_{\tilde{q_j}} \tau_{i_1,\cdots,i_{2k};0,\cdots,2k-1}^{\cdots,n_j-\tilde{q_j},\cdots}):\cdots]=[\cdots:
(\underset{\{l_1,\cdots,l_{2k}\}= \{0,\cdots,2k-1\}}{\sum} \pm\prod \bar{\nu}_{m_i-l_i})\cdot (\underset{\tilde{q}_j}{\sum}^{(2^k-2k)} \prod \bar{\nu}_{\tilde{q_j}} \bar{\tau}_{i_1,\cdots,i_{2k};0,\cdots,2k-1}^{\cdots,n_j-\tilde{q_j},\cdots}):\cdots]$, we have that $[\cdots:$
$\underset{\{l_1,\cdots,l_{2k}\}=
\{0,\cdots,2k-1\}}{\sum} \pm\prod \nu_{m_i-l_i}:\cdots]= [\cdots:\underset{\{l_1,\cdots,l_{2k}\}=\{0,\cdots,2k-1\}}{\sum} \pm\prod \bar{\nu}_{m_i-l_i}:\cdots]$. Since all the monomials of degree $2k$ are linear combinations of $\underset{\{l_1,\cdots,l_{2k}\}=\{0,\cdots,2k-1\}}{\sum} \pm\prod \nu_{m_i-l_i}$, the previous equation implies that $[\cdots: \nu_0^{i_0} \cdot \cdots \nu_{d-k}^{i_{d-k}}:\cdots]=[\cdots: \bar{\nu}_0^{i_0}\cdot \cdots \bar{\nu}_{d-k}^{i_{d-k}}:\cdots]$, which  is exactly the $2k$-uple embedding of $\P^{d-k}$ in $\P^{\binom{d+k}{2k}-1}$, so $[\nu_0:\cdots:\nu_{d-k}]=[\bar{\nu}_0:\cdots:\bar{\nu}_{d-k}]$.

Arrange $\tau_{i_1,\cdots,i_{2k};0,\cdots,2k-1}^{J_k}$ in the lexicographical order. Let $\tau_{i_1^{\circ},\cdots,i_{2k}^{\circ};0,\cdots,2k-1}^{\cdots,n_j^{\circ},\cdots}$ be the last non-zero element in the coordinate of $P$ among all $\tau_{i_1,\cdots,i_{2k};0,\cdots,2k-1}^{J_k}$. Let $\nu_q$ be the last non-zero element in the coordinate of $P$ among $\nu_0,\cdots,\nu_{d-k}$. Let $n_j=n_j^{\circ}+q$. Then $\underset{\tilde{q}_j}{\sum}^{(2^k-2k)} \prod \nu_{\tilde{q_j}} \tau_{i_1^{\circ},\cdots,i_{2k}^{\circ};0,\cdots,2k-1}^{\cdots,n_j-\tilde{q_j},\cdots}=\nu_q^{2^k-2k} \tau_{i_1^{\circ},\cdots,i_{2k}^{\circ};0,\cdots,2k-1}^{\cdots,n_j^{\circ},\cdots} \ne 0$. From $[\cdots:$

\noindent $(\underset{\{l_1,\cdots,l_{2k}\}=
\{0,\cdots,2k-1\}}{\sum} \pm\prod \nu_{m_i-l_i})\cdot (\underset{\tilde{q}_j}{\sum}^{(2^k-2k)} \prod \nu_{\tilde{q_j}} \tau_{i_1,\cdots,i_{2k};0,\cdots,2k-1}^{\cdots,n_j-\tilde{q_j},\cdots}):\cdots]=[\cdots:$

\noindent $(\underset{\{l_1,\cdots,l_{2k}\}= \{0,\cdots,2k-1\}}{\sum} \pm\prod \bar{\nu}_{m_i-l_i})\cdot (\underset{\tilde{q}_j}{\sum}^{(2^k-2k)} \prod \bar{\nu}_{\tilde{q_j}} \bar{\tau}_{i_1,\cdots,i_{2k};0,\cdots,2k-1}^{\cdots,n_j-\tilde{q_j},\cdots}):\cdots]$, we have that $[\cdots:$

\noindent $(\underset{\{l_1,\cdots,l_{2k}\}=
\{0,\cdots,2k-1\}}{\sum} \pm\prod \nu_{m_i-l_i})\cdot (\underset{\tilde{q}_j}{\sum}^{(2^k-2k)} \prod \nu_{\tilde{q_j}} \tau_{i_1^{\circ},\cdots,i_{2k}^{\circ};0,\cdots,2k-1}^{\cdots,n_j-\tilde{q_j},\cdots}):\cdots]=[\cdots:$

\noindent $(\underset{\{l_1,\cdots,l_{2k}\}= \{0,\cdots,2k-1\}}{\sum} \pm\prod \bar{\nu}_{m_i-l_i})\cdot (\underset{\tilde{q}_j}{\sum}^{(2^k-2k)} \prod \bar{\nu}_{\tilde{q_j}} \bar{\tau}_{i_1^{\circ},\cdots,i_{2k}^{\circ};0,\cdots,2k-1}^{\cdots,n_j-\tilde{q_j},\cdots}):\cdots]$, so $\underset{\tilde{q}_j}{\sum}^{(2^k-2k)} \prod \bar{\nu}_{\tilde{q_j}}$

\noindent $\bar{\tau}_{i_1^{\circ},\cdots,i_{2k}^{\circ};0,\cdots,2k-1}^{\cdots,n_j-\tilde{q_j},\cdots} \ne 0$. Since $[\nu_0:\cdots:\nu_{d-k}]=[\bar{\nu}_0:\cdots:\bar{\nu}_{d-k}]$, we have that $\bar{\nu}_q$ is the last non-zero element in the coordinate of $\bar{P}$ among $\bar{\nu}_0,\cdots,\bar{\nu}_{d-k}$. It follows from

\noindent $\underset{\tilde{q}_j}{\sum}^{(2^k-2k)} \prod \bar{\nu}_{\tilde{q_j}} $
$\bar{\tau}_{i_1^{\circ},\cdots,i_{2k}^{\circ};0,\cdots,2k-1}^{\cdots,n_j-\tilde{q_j},\cdots} \ne 0$ that the last non-zero element in the coordinate of $\bar{P}$ among all $\bar{\tau}_{i_1,\cdots,i_{2k};0,\cdots,2k-1}^{J_k}$ is equal to or strictly behind $\bar{\tau}_{i_1^{\circ},\cdots,i_{2k}^{\circ};0,\cdots,2k-1}^{\cdots,n_j^{\circ},\cdots}$. From the construction and the symmetric role of $\nu_q$ and $\bar{\nu}_q$, we see that the last non-zero element in the coordinate of $\bar{P}$ among all $\bar{\tau}_{i_1,\cdots,i_{2k};0,\cdots,2k-1}^{J_k}$ is $\bar{\tau}_{i_1^{\circ},\cdots,i_{2k}^{\circ};0,\cdots,2k-1}^{\cdots,n_j^{\circ},\cdots}$. So $\underset{\tilde{q}_j}{\sum}^{(2^k-2k)} \prod \bar{\nu}_{\tilde{q_j}} \bar{\tau}_{i_1^{\circ},\cdots,i_{2k}^{\circ};0,\cdots,2k-1}^{\cdots,n_j-\tilde{q_j},\cdots}$

\noindent $=\bar{\nu}_q^{2^k-2k} \bar{\tau}_{i_1^{\circ},\cdots,i_{2k}^{\circ};0,\cdots,2k-1}^{\cdots,n_j^{\circ},\cdots} \ne 0$. From $[\cdots:(\underset{\{l_1,\cdots,l_{2k}\}=
\{0,\cdots,2k-1\}}{\sum} \pm\prod \nu_{m_i-l_i})\cdot (\underset{\tilde{q}_j}{\sum}^{(2^k-2k)} \prod \nu_{\tilde{q_j}}$

\noindent $\tau_{i_1,\cdots,i_{2k};0,\cdots,2k-1}^{\cdots,n_j-\tilde{q_j},\cdots}):\cdots]=[\cdots:
(\underset{\{l_1,\cdots,l_{2k}\}= \{0,\cdots,2k-1\}}{\sum} \pm\prod \bar{\nu}_{m_i-l_i})\cdot (\underset{\tilde{q}_j}{\sum}^{(2^k-2k)} \prod \bar{\nu}_{\tilde{q_j}}$

\noindent $\bar{\tau}_{i_1,\cdots,i_{2k};0,\cdots,2k-1}^{\cdots,n_j-\tilde{q_j},\cdots}):\cdots]$, we have 
$[\underset{\tilde{q}_j}{\sum}^{(2^k-2k)}$
$\prod \nu_{\tilde{q_j}} \tau_{i_1,\cdots,i_{2k};0,\cdots,2k-1}^{\cdots,n_j-\tilde{q_j},\cdots}: \nu_q^{2^k-2k} \tau_{i_1^{\circ},\cdots,i_{2k}^{\circ};0,\cdots,2k-1}^{\cdots,n_j^{\circ},\cdots} ]$

\noindent $=[ \underset{\tilde{q}_j}{\sum}^{(2^k-2k)} \prod \bar{\nu}_{\tilde{q_j}} \bar{\tau}_{i_1,\cdots,i_{2k};0,\cdots,2k-1}^{\cdots,n_j-\tilde{q_j},\cdots} : \bar{\nu}_q^{2^k-2k}$
$\bar{\tau}_{i_1^{\circ},\cdots,i_{2k}^{\circ};0,\cdots,2k-1}^{\cdots,n_j^{\circ},\cdots} ]$. We will use this relation and $[\nu_0:\cdots:\nu_{d-k}]=[\bar{\nu}_0:\cdots:\bar{\nu}_{d-k}]$ to show that $[\cdots:\tau_{i_1,\cdots,i_{2k};0,\cdots,2k-1}^{J_k}:\cdots]= [\cdots:\bar{\tau}_{i_1,\cdots,i_{2k};0,\cdots,2k-1}^{J_k}:\cdots]$. For simplicity, let $\nu_q=\bar{\nu}_q=
\tau_{i_1^{\circ},\cdots,i_{2k}^{\circ};0,\cdots,2k-1}^{\cdots,n_j^{\circ},\cdots}= \bar{\tau}_{i_1^{\circ},\cdots,i_{2k}^{\circ};0,\cdots,2k-1}^{\cdots,n_j^{\circ},\cdots}=1$. Then $\nu_i=\bar{\nu}_i (i=0,\cdots,d-k)$ and $\underset{\tilde{q}_j}{\sum}^{(2^k-2k)} \prod \nu_{\tilde{q_j}} \tau_{i_1,\cdots,i_{2k};0,\cdots,2k-1}^{\cdots,n_j-\tilde{q_j},\cdots}=\underset{\tilde{q}_j}{\sum}^{(2^k-2k)} \prod \bar{\nu}_{\tilde{q_j}}$
$\bar{\tau}_{i_1,\cdots,i_{2k};0,\cdots,2k-1}^{\cdots,n_j-\tilde{q_j},\cdots}$. Now let us show that $\tau_{i_1,\cdots,i_{2k};0,\cdots,2k-1}^{J_k}=\bar{\tau}_{i_1,\cdots,i_{2k};0,\cdots,2k-1}^{J_k}$.

Let $\tau_{i'_1,\cdots,i'_{2k};0,\cdots,2k-1}^{J'_k}$ and $\bar{\tau}_{i'_1,\cdots,i'_{2k};0,\cdots,2k-1}^{J'_k}$ be the coordinate component of $P$ and $\bar{P}$  immediately in front of $\tau_{i_1^{\circ},\cdots,i_{2k}^{\circ};0,\cdots,2k-1}^{\cdots,n_j^{\circ},\cdots}$ and $\bar{\tau}_{i_1^{\circ},\cdots,i_{2k}^{\circ};0,\cdots,2k-1}^{\cdots,n_j^{\circ},\cdots}$ in $\P^{\binom{r+1}{2}^{k} \cdot \overset{k-2}{\underset{j=1}
{\prod}}\binom{2k-1-j}{2(k-1-j)}^{2^{j-1}} 
-1}$, respectively. If $(i'_1,\cdots,i'_{2k})< (i_1^{\circ},\cdots,i_{2k}^{\circ})$, then $J'_k=(\cdots, n_j^{\circ}, \cdots)$. Pick $n_j=q+ n_j^{\circ}$. Then $\underset{\tilde{q}_j}{\sum}^{(2^k-2k)} \prod \nu_{\tilde{q_j}} \tau_{i'_1,\cdots,i'_{2k};0,\cdots,2k-1}^{\cdots,n_j-\tilde{q_j},\cdots}= \nu_q^{2^k-2k} \tau_{i'_1,\cdots,i'_{2k};0,\cdots,2k-1}^{\cdots,n_j^{\circ},\cdots}= \tau_{i'_1,\cdots,i'_{2k};0,\cdots,2k-1}^{\cdots,n_j^{\circ},\cdots}$ because any $\tilde{q_j}< q$ will give an element between $\tau_{i'_1,\cdots,i'_{2k};0,\cdots,2k-1}^{J'_k}$ and $\tau_{i_1^{\circ},\cdots,i_{2k}^{\circ};0,\cdots,2k-1}^{\cdots,n_j^{\circ},\cdots}$, which is impossible due to the choice of $\tau_{i'_1,\cdots,i'_{2k};0,\cdots,2k-1}^{J'_k}$. Similarly for $n_j=q+ n_j^{\circ}$, we have $\underset{\tilde{q}_j}{\sum}^{(2^k-2k)} \prod \bar{\nu}_{\tilde{q_j}} \bar{\tau}_{i'_1,\cdots,i'_{2k};0,\cdots,2k-1}^{\cdots,n_j-\tilde{q_j},\cdots}= \bar{\nu}_q^{2^k-2k} \bar{\tau}_{i'_1,\cdots,i'_{2k};0,\cdots,2k-1}^{\cdots,n_j^{\circ},\cdots}=\bar{\tau}_{i'_1,\cdots,i'_{2k};0,\cdots,2k-1}^{\cdots,n_j^{\circ},\cdots}$. Now $\underset{\tilde{q}_j}{\sum}^{(2^k-2k)} \prod \nu_{\tilde{q_j}} \tau_{i'_1,\cdots,i'_{2k};0,\cdots,2k-1}^{\cdots,n_j-\tilde{q_j},\cdots}= \underset{\tilde{q}_j}{\sum}^{(2^k-2k)} \prod \bar{\nu}_{\tilde{q_j}} \bar{\tau}_{i'_1,\cdots,i'_{2k};0,\cdots,2k-1}^{\cdots,n_j-\tilde{q_j},\cdots}$ gives
$\tau_{i'_1,\cdots,i'_{2k};0,\cdots,2k-1}^{\cdots,n_j^{\circ},\cdots}= \bar{\tau}_{i'_1,\cdots,i'_{2k};0,\cdots,2k-1}^{\cdots,n_j^{\circ},\cdots}$. If $(i'_1,\cdots,i'_{2k})= (i_1^{\circ},\cdots,i_{2k}^{\circ})$, then $J'_k =(\cdots, n'_j, \cdots)< (\cdots, n_j^{\circ}, \cdots)$. Pick $n_j=n'_j+q$. Then $\underset{\tilde{q}_j}{\sum}^{(2^k-2k)} \prod \nu_{\tilde{q_j}} \tau_{i_1^{\circ},\cdots,i_{2k}^{\circ};0,\cdots,2k-1}^{\cdots,n_j-\tilde{q_j},\cdots}= \nu_q^{2^k-2k} \tau_{i_1^{\circ},\cdots,i_{2k}^{\circ};0,\cdots,2k-1}^{\cdots,n'_j,\cdots}+ \underset{\tilde{q}_j}{\sum'}^{(2^k-2k)}$

\noindent $\prod \nu_{\tilde{q_j}} \tau_{i_1^{\circ},\cdots,i_{2k}^{\circ};0,\cdots,2k-1}^{\cdots,n_j-\tilde{q_j},\cdots}$, where $\sum'$ means that the sum is over all $\tilde{q}_j$ with at least one less than $q$. Each term $\tau_{i_1^{\circ},\cdots,i_{2k}^{\circ};0,\cdots,2k-1}^{\cdots,n_j-\tilde{q_j},\cdots}$ in $\sum'$ is behind $\tau_{i_1^{\circ},\cdots,i_{2k}^{\circ};0,\cdots,2k-1}^{\cdots,n'_j,\cdots}$, so it is either $\tau_{i_1^{\circ},\cdots,i_{2k}^{\circ};0,\cdots,2k-1}^{\cdots,n_j^{\circ},\cdots}$ or zero and it agrees with $\bar{\tau}_{i_1^{\circ},\cdots,i_{2k}^{\circ};0,\cdots,2k-1}^{\cdots,n_j-\tilde{q_j},\cdots}$, thus $\underset{\tilde{q}_j}{\sum'}^{(2^k-2k)} \prod \nu_{\tilde{q_j}} \tau_{i_1^{\circ},\cdots,i_{2k}^{\circ};0,\cdots,2k-1}^{\cdots,n_j-\tilde{q_j},\cdots}= \underset{\tilde{q}_j}{\sum'}^{(2^k-2k)} $

\noindent $\prod \bar{\nu}_{\tilde{q_j}} \bar{\tau}_{i_1^{\circ},\cdots,i_{2k}^{\circ};0,\cdots,2k-1}^{\cdots,n_j-\tilde{q_j},\cdots}$. Hence 
$\tau_{i_1^{\circ},\cdots,i_{2k}^{\circ};0,\cdots,2k-1}^{\cdots,n'_j,\cdots}=\nu_q^{2^k-2k} \tau_{i_1^{\circ},\cdots,i_{2k}^{\circ};0,\cdots,2k-1}^{\cdots,n'_j,\cdots}= \underset{\tilde{q}_j}{\sum}^{(2^k-2k)} $

\noindent $\prod \nu_{\tilde{q_j}} \tau_{i_1^{\circ},\cdots,i_{2k}^{\circ};0,\cdots,2k-1}^{\cdots,n_j-\tilde{q_j},\cdots}- \underset{\tilde{q}_j}{\sum'}^{(2^k-2k)} \prod \nu_{\tilde{q_j}} \tau_{i_1^{\circ},\cdots,i_{2k}^{\circ};0,\cdots,2k-1}^{\cdots,n_j-\tilde{q_j},\cdots}$
$= \underset{\tilde{q}_j}{\sum}^{(2^k-2k)} \prod \bar{\nu}_{\tilde{q_j}} \bar{\tau}_{i_1^{\circ},\cdots,i_{2k}^{\circ};0,\cdots,2k-1}^{\cdots,n_j-\tilde{q_j},\cdots}- $

\noindent $\underset{\tilde{q}_j}{\sum'}^{(2^k-2k)} \prod \bar{\nu}_{\tilde{q_j}} \bar{\tau}_{i_1^{\circ},\cdots,i_{2k}^{\circ};0,\cdots,2k-1}^{\cdots,n_j-\tilde{q_j},\cdots}=\bar{\nu}_q^{2^k-2k} \bar{\tau}_{i_1^{\circ},\cdots,i_{2k}^{\circ};0,\cdots,2k-1}^{\cdots,n'_j,\cdots}=$
$\bar{\tau}_{i_1^{\circ},\cdots,i_{2k}^{\circ};0,\cdots,2k-1}^{\cdots,n'_j,\cdots}$. Therefore in either case we have $\tau_{i'_1,\cdots,i'_{2k};0,\cdots,2k-1}^{J'_k}=\bar{\tau}_{i'_1,\cdots,i'_{2k};0,\cdots,2k-1}^{J'_k}$.

Suppose that the correspondening coordinate components of $P$ and $\bar{P}$ behind

\noindent $\tau_{i_1,\cdots,i_{2k};0,\cdots,2k-1}^{\cdots,n_j,\cdots}$ and $\bar{\tau}_{i_1,\cdots,i_{2k};0,\cdots,2k-1}^{\cdots,n_j,\cdots}$ in $\P^{\binom{r+1}{2}^{k} \cdot \overset{k-2}{\underset{j=1}
{\prod}}\binom{2k-1-j}{2(k-1-j)}^{2^{j-1}} 
-1}$ agree, we will show that $\tau_{i_1,\cdots,i_{2k};0,\cdots,2k-1}^{\cdots,n_j,\cdots}=\bar{\tau}_{i_1,\cdots,i_{2k};0,\cdots,2k-1}^{\cdots,n_j,\cdots}$ as well.

Replacing $n_j$ with $n_j+q$ in the formula $\underset{\tilde{q}_j}{\sum}^{(2^k-2k)} \prod \nu_{\tilde{q_j}} \tau_{i_1,\cdots,i_{2k};0,\cdots,2k-1}^{\cdots,n_j-\tilde{q_j},\cdots}=\underset{\tilde{q}_j}{\sum}^{(2^k-2k)} \prod \bar{\nu}_{\tilde{q_j}}$

\noindent $\bar{\tau}_{i_1,\cdots,i_{2k};0,\cdots,2k-1}^{\cdots,n_j-\tilde{q_j},\cdots}$, we have that $ \underset{\tilde{q}_j}{\sum}^{(2^k-2k)} \prod \nu_{\tilde{q_j}} \tau_{i_1,\cdots,i_{2k};0,\cdots,2k-1}^{\cdots,n_j+q-\tilde{q_j},\cdots}= \nu_q^{2^k-2k} \tau_{i_1,\cdots,i_{2k};0,\cdots,2k-1}^{\cdots,n_j,\cdots}$

\noindent $+ \underset{\tilde{q}_j}{\sum'}^{(2^k-2k)} \prod \nu_{\tilde{q_j}} \tau_{i_1,\cdots,i_{2k};0,\cdots,2k-1}^{\cdots,n_j+q-\tilde{q_j},\cdots}= \underset{\tilde{q}_j}{\sum}^{(2^k-2k)} \prod \bar{\nu}_{\tilde{q_j}} \bar{\tau}_{i_1,\cdots,i_{2k};0,\cdots,2k-1}^{\cdots,n_j+q-\tilde{q_j},\cdots}=\bar{\nu}_q^{2^k-2k} \bar{\tau}_{i_1,\cdots,i_{2k};0,\cdots,2k-1}^{\cdots,n_j,\cdots}$

\noindent $+ \underset{\tilde{q}_j}{\sum'}^{(2^k-2k)} \prod \bar{\nu}_{\tilde{q_j}} \bar{\tau}_{i_1,\cdots,i_{2k};0,\cdots,2k-1}^{\cdots,n_j+q-\tilde{q_j},\cdots}$, where $\sum'$ again represents the summation over all the $\tilde{q_j}$ with at least one less than $q$. By our inductive assumption, $\underset{\tilde{q}_j}{\sum'}^{(2^k-2k)} \prod \nu_{\tilde{q_j}} \tau_{i_1,\cdots,i_{2k};0,\cdots,2k-1}^{\cdots,n_j+q-\tilde{q_j},\cdots}$

\noindent $= \underset{\tilde{q}_j}{\sum'}^{(2^k-2k)} \prod \bar{\nu}_{\tilde{q_j}} \bar{\tau}_{i_1,\cdots,i_{2k};0,\cdots,2k-1}^{\cdots,n_j+q-\tilde{q_j},\cdots}$. So   
the previous equation implies that $\tau_{i_1,\cdots,i_{2k};0,\cdots,2k-1}^{\cdots,n_j,\cdots}$

\noindent $=\nu_q^{2^k-2k} \tau_{i_1,\cdots,i_{2k};0,\cdots,2k-1}^{\cdots,n_j,\cdots}= \bar{\nu}_q^{2^k-2k} \bar{\tau}_{i_1,\cdots,i_{2k};0,\cdots,2k-1}^{\cdots,n_j,\cdots}=\bar{\tau}_{i_1,\cdots,i_{2k};0,\cdots,2k-1}^{\cdots,n_j,\cdots}$. This completes the proof that the coordinate components of $P$ and $\bar{P}$  agree in $\P^{\binom{r+1}{2}^{k} \cdot \overset{k-2}{\underset{j=1}{\prod}}\binom{2k-1-j}{2(k-1-j)}^{2^{j-1}} 
-1}$.    

Essentially in the same way we can show that the coordinate components of $P$ and $\bar{P}$ agree in 
$\P^{\binom{r+1}{2}^{s} \cdot \binom{k+s}{2s}\cdot \overset{k-2}{\underset{j=1}
{\prod}}\binom{k+(s-1-j)}{2(s-1-j)}^{2^{j-1}} 
-1}$ for $s=3,\cdots,k-1$, $\P^{\binom{r+1}{2} \cdot \binom{k+1}{2}-1}$, $\P^{\binom{r+1}{2} \cdot \binom{k+1}{2}-1}$, and $\P^{\binom{r+1}{2}^2 \cdot \binom{k+2}{4}-1}$, respectively. We leave it as an exercise to the interest readers. From the definition of $\Phi_{k+1}: \map \P^{(k+1)(r+1)-1} \times \P^{d-k}. R_{k+1} \subset \P^{(d+1)(r+1)-1}.$ and $\nu_i=\bar{\nu}_i ( i=0,\cdots,d-k)$, it is easy to see that the coordinate components of $P$ and $\bar{P}$ in $\P^{(k+1)(r+1)-1}$ agree as well. This completes the proof that $F$ is injective.

Since $F$ is an injective morphism between projective varieties, it is finite and hence closed. That $\pi^{-1} (R_{k+1} \setminus R_k) \subset F(\Gamma_{\varphi'_{k}} \times \P^{d-k})$ and $F$ closed implies $\widetilde{R}_{k+1}$, the closure of $\pi^{-1} (R_{k+1} \setminus R_k)$ in $\Gamma_{\varphi_{k}}$, is a subset of $F(\Gamma_{\varphi'_{k}} \times \P^{d-k})$, so $F$ is onto. 

Now we have proved that $F: \map \Gamma_{\varphi'_{k}} \times \P^{d-k}. \widetilde{R}_{k+1}.$ is a bijection. In order to show that $F$ is actually an isomorphism, it is enough to show that $F$ is a closed immersion. 

By a local criterion on closed immersion (see, for example, Proposition 7.3, Chapter II in Hartshorne $[3]$), we need to verify that the coordinate functions in $F$ separate points and tangent vectors. That the coordinate functions in $F$ separate points has been shown above. Now let us show that they also separate tangent vectors.

Since $F$ is a homeomorphism onto $\widetilde{R}_{k+1}$, to show that the coordinate functions in $F$ separate tangent vectors, we only need to show that the morphsim of sheaves $\map \CO_{\widetilde{R}_{k+1}}. F_{*} \CO_{\Gamma_{\varphi'_k} \times \P^{d-k}}.$ is surjective. We check this surjectivity on stalks.

From the proof of injectivity, we see that $F$ restricted to $\P^{d-k} \times \P^{\binom{r+1}{2}^{k} \cdot \overset{k-2}{\underset{j=1}
{\prod}}\binom{2k-1-j}{2(k-1-j)}^{2^{j-1}} 
-1}$ gives a morphism $G: \map \P^{d-k} \times \P^{\binom{r+1}{2}^{k} \cdot \overset{k-2}{\underset{j=1}
{\prod}}\binom{2k-1-j}{2(k-1-j)}^{2^{j-1}} 
-1}. \P^{\binom{r+1}{2}^{k} \cdot \binom{d+k}{2k}\cdot \overset{k-2}{\underset{j=1}
{\prod}}\binom{d+(k-1-j)}{2(k-1-j)}^{2^{j-1}} 
-1}.$, where $G([\nu_0:\cdots:\nu_{d-k}];[\cdots:\tau_{i_1,\cdots,i_{2k};0,\cdots,2k-1}^{\cdots,n_j,\cdots}:\cdots])=[\cdots:(\underset{\{l_1,\cdots,l_{2k}\}=
\{0,\cdots,2k-1\}}{\sum} \pm\prod \nu_{m_i-l_i})\cdot (\underset{\tilde{q}_j}{\sum}^{(2^k-2k)} \prod \nu_{\tilde{q_j}}$
$\tau_{i_1,\cdots,i_{2k};0,\cdots,2k-1}^{\cdots,n_j-\tilde{q_j},\cdots}):\cdots]$. Let $\tilde{G}: \P^{\binom{d+k}{2k}-1} \times $

\noindent $\P^{\binom{r+1}{2}^{k} \cdot \overset{k-2}{\underset{j=1}
{\prod}}\binom{d+(k-1-j)}{2(k-1-j)}^{2^{j-1}} 
-1} \longrightarrow \P^{\binom{r+1}{2}^{k} \cdot \binom{d+k}{2k}\cdot \overset{k-2}{\underset{j=1}
{\prod}}\binom{d+(k-1-j)}{2(k-1-j)}^{2^{j-1}} 
-1}$ be the standard Serge embedding. Then $G$ can be decomposed as $\tilde{G} \circ H$, where $H: \P^{d-k} \times \P^{\binom{r+1}{2}^{k} \cdot \overset{k-2}{\underset{j=1}
{\prod}}\binom{2k-1-j}{2(k-1-j)}^{2^{j-1}} 
-1} \longrightarrow \P^{\binom{d+k}{2k}-1} \times $
$\P^{\binom{r+1}{2}^{k} \cdot \overset{k-2}{\underset{j=1}
{\prod}}\binom{2k-1-j}{2(k-1-j)}^{2^{j-1}} 
-1}$ is the morphism sending $[\nu_0:\cdots:\nu_{d-k}]; [\cdots: $

\noindent $\tau_{i_1,\cdots,i_{2k};0,\cdots,2k-1}^{\cdots,n_j,\cdots}:\cdots]$ into $[\cdots:\underset{\{l_1,\cdots,l_{2k}\}=
\{0,\cdots,2k-1\}}{\sum} \pm\prod \nu_{m_i-l_i}:\cdots]; [\cdots:\underset{\tilde{q}_j}{\sum}^{(2^k-2k)} \prod \nu_{\tilde{q_j}}$
$\tau_{i_1,\cdots,i_{2k};0,\cdots,2k-1}^{\cdots,n_j-\tilde{q_j},\cdots}:\cdots]$. Since $\tilde{G}$ is an embedding when restricted to the image of $H$, then $\underset{\{l_1,\cdots,l_{2k}\}=
\{0,\cdots,2k-1\}}{\sum} \pm\prod \nu_{m_i-l_i}$ and $\underset{\tilde{q}_j}{\sum}^{(2^k-2k)} \prod \nu_{\tilde{q_j}} \tau_{i_1,\cdots,i_{2k};0,\cdots,2k-1}^{\cdots,n_j-\tilde{q_j},\cdots}$ are locally regular functions in terms of 
$u_{i_1,\cdots,i_{2k};m_1,\cdots,m_{2k}}^{\cdots,n_j,\cdots}$.
We have showed that all the monomials of degree $2k$ in $\nu_0,\cdots,\nu_{d-k}$ are linear combination of $\underset{\{l_1,\cdots,l_{2k}\}=
\{0,\cdots,2k-1\}}{\sum} \pm\prod \nu_{m_i-l_i}$. So locally $\nu_0,\cdots,\nu_{d-k}$ are regular functions of $\underset{\{l_1,\cdots,l_{2k}\}=
\{0,\cdots,2k-1\}}{\sum} \pm\prod \nu_{m_i-l_i}$. Therefore $\nu_0,\cdots,\nu_{d-k}$ are locally regular functions in terms of $u_{i_1,\cdots,i_{2k};m_1,\cdots,m_{2k}}^{\cdots,n_j,\cdots}$.

Given any point $Q$ in $\widetilde{R}_{k+1}$, let $F^{-1} (Q)=P$. We need to show that  
$\map \CO_{\widetilde{R}_{k+1}, Q}. \CO_{\Gamma_{\varphi'_k} \times \P^{d-k}, P}.$ is surjective.

Pick affine neighborhoods around $P$ and $Q$. For convenience we still use the same letters to represent coordinates but the reader should be aware of that they are affine now. Consider the following linear equations in unknown variables $\mu_{i,j},\tau_{i_1,i_2;m_1,m_2},$

\noindent $\tau_{i_1,\cdots,i_4;m_1,\cdots,m_4}$ and $\tau_{i_1,\cdots,i_{2s};m_1,\cdots,m_{2s}}^{J_s}, s=3,\cdots,k$:

\begin{equation}
\begin{split}
\underset{p+q=j}{\sum} \nu_q \mu_{ip}=s_{ij}\cr
\underset{p_1+q_1=m_1}{\sum} \underset{p_2+q_2=m_2}{\sum} \nu_{q_1} \nu_{q_2} \tau_{i_1i_2,p_1p_2}=u_{i_1,i_2;m_1,m_2}\cr
\overset{m_1}{\underset{q_1=0}{\sum}} \overset{m_2}{\underset{q_2=0}{\sum}}  \overset{m_3}{\underset{q_3=0}{\sum}} \overset{m_4}{\underset{q_4=0}{\sum}} \nu_{q_1} 
\nu_{q_2}\nu_{q_3} \nu_{q_4} \tau_{i_1i_2i_3i_4;m_1-q_1,m_2-q_2,m_3-q_3,m_4-q_4}=u_{i_1,\cdots,i_4;m_1,\cdots,m_4}\cr 
\underset{q_i}{\sum}^{(2s)} \underset{\tilde{q}_j}{\sum}^{(2^s-2s)} \prod \nu_{q_i} \cdot \prod \nu_{\tilde{q_j}} \tau_{i_1,\cdots,i_{2s};m_1-q_1,\cdots,m_{2s}-q_{2s}}^{\cdots,n_j-\tilde{q_j},\cdots}=u_{i_1,\cdots,i_{2s};m_1,\cdots,m_{2s}}^{\cdots,n_j,\cdots}\cr
s=3,\cdots,k\cr 
\end{split}
\end{equation}

By the injectivity of $F$, for part of the coordinates $\nu_0,\cdots,\nu_{d-k}$ in $P$ and the coordinates $s_{ij},u_{i_1,i_2;m_1,m_2}, u_{i_1,\cdots,i_4;m_1,\cdots,m_4}, u_{i_1,\cdots,i_{2s};m_1,\cdots,m_{2s}}^{\cdots,n_j,\cdots},s=3,\cdots,k$ in $Q$, the above linear equations have a unique solution, which implies that the coefficient matrix is of full column rank and $\mu_{i,j},\tau_{i_1,i_2;m_1,m_2},$
$\tau_{i_1,\cdots,i_4;m_1,\cdots,m_4}$ and $\tau_{i_1,\cdots,i_{2s};m_1,\cdots,m_{2s}}^{J_s}, s=3,\cdots,k$ are given by quotient of determinants involving $\nu_i,i=0,\cdots,$
$d-k$ and $s_{ij},u_{i_1,i_2;m_1,m_2}, u_{i_1,\cdots,i_4;m_1,\cdots,m_4},$

\noindent $u_{i_1,\cdots,i_{2s};m_1,\cdots,m_{2s}}^{\cdots,n_j,\cdots},s=3,\cdots,k$. Since full column rank is an open condition, by shrinking the affine neighborhood of $P$ if necessary, we can assume that the coefficient matrix is of full column rank everywhere in the neighborhood. Hence $\mu_{i,j},\tau_{i_1,i_2;m_1,m_2},\tau_{i_1,\cdots,i_4;m_1,\cdots,m_4}$ and $\tau_{i_1,\cdots,i_{2s};m_1,\cdots,m_{2s}}^{J_s}, s=3,\cdots,k$ can be given by quotient of determinants involving variables $\nu_i,i=0,\cdots,d-k$,  $s_{ij}$,$u_{i_1,i_2;m_1,m_2}$, $u_{i_1,\cdots,i_4;m_1,\cdots,m_4},$ and $u_{i_1,\cdots,i_{2s};m_1,\cdots,m_{2s}}^{\cdots,n_j,\cdots},$ $s=3,\cdots$$,k$ locally. We have known that locally $\nu_i,i=0,\cdots,d-k$  are regular functions 
in terms of $u_{i_1,\cdots,i_{2k};m_1,\cdots,m_{2k}}^{\cdots,n_j,\cdots}$, so locally $\mu_{i,j},\tau_{i_1,i_2;m_1,m_2},\tau_{i_1,\cdots,i_4;m_1,\cdots,m_4}$ and $\tau_{i_1,\cdots,i_{2s};m_1,\cdots,m_{2s}}^{J_s},$
$s=3,\cdots$ $,k$ are regular functions in terms of $s_{ij},u_{i_1,i_2;m_1,m_2}, u_{i_1,\cdots,i_4;m_1,\cdots,m_4}$, and $u_{i_1,\cdots,i_{2s};m_1,\cdots,m_{2s}}^{\cdots,n_j,\cdots},$
$s=3$,$\cdots,k$. This completes the proof that $\map \CO_{\widetilde{R}_{k+1}}. F_{*} \CO_{\Gamma_{\varphi'_k} \times \P^{d-k}}.$ is surjective. So $F$ is an immersion.

Since $F: \map \Gamma_{\varphi'_{k}} \times \P^{d-k}. \widetilde{R}_{k+1}.$ is a bijection and immersion, it is an isomorphism. The smoothness of $\widetilde{R}_{k+1}$ follows from that of $\Gamma_{\varphi'_{k}} \times \P^{d-k}$. Therefore the $(k+1)$-th blow-up can be carried out. This completes the proof that the iterated blow-ups along the strata (or its proper transformations) of $M_{d}(\P^r) \setminus M_{d}^{\circ}(\P^r)$ can be done. The final outcome of the iterated process is $\Gamma_{\varphi_d}$, whose boundary consists of an exceptional divisor $E_d$ and $d-1$ proper transformation of the exceptional divisors $\tilde{E}_i, i=1, \cdots, d-1$ appearing in the first $(d-1)$ blow-ups. To complete our proof of theorem $1.1$, we only need to show that the union of $E_d$ and $\tilde{E}_i, i=1, \cdots, d-1$ are normal crossing divisors.

For $d=1$, we have known that the boundary consists of an exceptional divisor, which is of course normal crossing. Suppose that the iterated blow-ups can be carried out for the space of holomorphic maps of degree less than $d$ and the final outcome has normal crossing boundary. From the proof of the smoothness of $\widetilde{R}_{k+1}$, we see that $F$ induces an isomorphism between $(E'_k \cap \underset{i \in I}{\cup} \tilde{E'}_i) \times \P^{d-k}$ and $\widetilde{R}_{k+1} \cap (E_k \cap \underset{i \in I}{\cup} \tilde{E}_i)$ for any subset $I$ of $\{1,\cdots,k-1\}$. So $\widetilde{R}_{k+1} \cap (E_k \cap \underset{i \in I}{\cup} \tilde{E}_i)$ is smooth, which implies that $\tilde{E}_{k+1}$ and $\tilde{E}_{i}, i=1, \cdots,k$ are normal crossing. Carrying out the whole iterated blow-ups will give that $E_d$ and $\tilde{E}_i, i=1, \cdots, d-1$ are normal crossing divisors.

\bigskip
\begin{center}  
{\bf References}
\end{center}
\bigskip

\begin{enumerate}
\item W. Fulton and R. MacPherson, A compactification of configuration spaces. Ann. of Math. (2) {\bf 139} (1994), no. 1, 183--225.

\item J. Harris, Algebraic Geometry. A first course. Corrected reprint of the 1992 original. Graduate Texts in Mathematics, no. {\bf 133}. Springer-Verlag, New York, 1995. xx+328 pp. 

\item R. Hartshorne, Algebraic Geometry. Graduate Texts in Mathematics, no. {\bf 52}. Springer-Verlag, New York-Heidelberg, 1977. xvi+496 pp.

\item Y. Hu, A compactification of open varieties. Trans. Amer. Math. Soc. 355 (2003), no. 12, 4737--4753.

\item K. Kaki\'e, The resultant of several homogeneous polynomials in two indeterminates. Proc. Amer. Math. Soc. {\bf 54} (1976), 1--7.

\item R. MacPherson and C. Procesi, Making conical compactifications wonderful. Selecta Math. (N.S.) {\bf 4} (1998), no. {\bf 1}, 125--139.

\item A. Ulyanov, Polydiagonal compactification of configuration spaces. J. Algebraic Geom. {\bf 11} (2002), no. {\bf 1}, 129--159.

\end{enumerate}
\medskip
\noindent Department of Mathematics, University of Georgia, Athens, GA 30602, USA.

\noindent jiayuan@uga.edu
\end{document}